\documentclass[11pt]{amsart}

\usepackage{amsmath,amsfonts,amssymb,amscd}
\begin{document}
\renewcommand{\thefootnote}{\fnsymbol{footnote}}
\pagestyle{plain}

%%%%%%%%%%%%%%%%%%%%%%%%
\title{{\small \text{Existence problem of extremal K\"ahler metrics}}}
\author{Toshiki Mabuchi${}^*$}
%\address{Department of Mathematics, Graduate School of Science, Osaka
%University,  Toyonaka, Osaka, 560-0043 Japan
%}
\maketitle
%$$
%\align
%&x=y
%\tag"(1.1)" \\
%&z=w \tag"(1.2)"
%\endalign
%$$
%%%%%%%%%%%%%%%%%%%%%%%%
\footnotetext{
%2010 {\it Mathematics Subject Classification.}
%Primary~32Q26; Secondary~14L24, 53C25.\\
${}^{*}$Supported 
by JSPS Grant-in-Aid for Scientific Research (A) No. 20244005.}
%%%%%%%%%%%%%%%%%%%%%%%%
\abstract
In this paper, we shall give some affirmative answer to an extremal K\"ahler version of 
the Yau-Tian-Donaldson Conjecture.
For a polarized algebraic manifold $(X,L)$, we choose a maximal algebraic torus  $T$ in the group $\operatorname{Aut}(X)$ of holomorphic automorphisms of $X$.
Then by the method as in \cite{M1}, the polarization class $c_1(L)$ will be shown to admit an extremal K\"ahler metric if $(X,L)$ is strongly K-stable relative to $T$.
\endabstract
%%%%%%%%%%%%%%%%%%%%%%

\section{Introduction}

In this paper, we discuss  the existence problem of extremal K\"ahler metrics on
 a {\it polarized algebraic manifold} $(X,L)$, i.e.,
a non-singular irreducible projective variety $X$, 
defined over $\Bbb C$, with a very ample line bundle $L$
on $X$.  Throughout this paper, we fix a pair $(X,L)$ as above 
with $n:= \dim^{}_{\Bbb C} X$, and assume that $T$ is a maximal algebraic torus in
$\operatorname{Aut}(M)$. However, the arguments below are valid also for an arbitrary algebraic torus 
sitting in $\operatorname{Aut}(M)$. 
The main purpose of this paper is to give an answer to some extremal K\"ahler version 
of  the Yau-Tian-Donaldson Conjecture:

\medskip\noindent
{\bf Main Theorem.} 
{\em For a polarized algebraic manifold $(X,L)$, the class $c_1(L)_{\Bbb R}$
 admits an extremal K\"ahler metric if $X$ is strongly K-stable relative to $T$.}

\section{Strong relative K-stability}

For the complex affine line 
 $ \Bbb A^1 := \{z\in \Bbb C\}$,
we consider the action of the algebraic torus $\Bbb G_m = \{\, t\in \Bbb C^*\}$ on $\Bbb A^1$ by multiplication of complex numbers
$$
\Bbb G_m\times \Bbb A^1 \to \Bbb A^1,
\qquad (t, z) \mapsto
t z.
$$
Fix a Hermitian metric $h$ for $L$ such that $\omega := c_1(L;h)$ is K\"ahler 
and that the maximal compact subgroup $T_c$ of $T$ acts isometrically on $(X, \omega )$. 
We then endow  
$V_{\ell}:= H^0(X,L^{\otimes\ell})$   
with the Hermitian metric $\rho_{\ell}$ defined by
$$
\langle \tau',\tau''\rangle_{\rho_{\ell}} \; :=\; \int_X \; (\tau', \tau'')_h \, \omega^n,
\qquad \tau',\tau'' \in V_{\ell},
\leqno{(2.1)}
$$
where $(\tau', \tau'')_h$ denotes the pointwise Hermitian inner product of $\tau'$ and $\tau''$ 
by the $\ell$-multiple of $h$.
 Let $ \Phi_{\ell} \,:\, X \, \hookrightarrow \, \Bbb P^*(V_{\ell})$
be the Kodaira embedding of $X$ associated to the complete linear system $|L^{\otimes\ell}|$ on $X$.
Let $X_{\ell}$ be its image $\Phi_{\ell}(X)$ in $\Bbb P^*(V_{\ell})$.
We then consider an algebraic group homomorphism
$$
\psi \, : \,\Bbb G_m \,\to \,\operatorname{SL}(V_{\ell}).
$$
such that the compact subgroup 
$S^1 \subset \Bbb G_m$ acts isometrically on $(V_{\ell}, \rho_{\ell})$.
Let $\mathcal{X}^{\psi}$ be the irreducible algebraic subvariety of $\Bbb A^1 \times \Bbb P^* (V_{\ell})$ obtained as 
the closure of $\cup_{z\in \Bbb C^*} \mathcal{X}^{\psi}_z$ in $\Bbb A^1 \times \Bbb P^* (V_{\ell})$ 
by setting
$$
\mathcal{X}^{\psi}_z := \{z\}\times\psi (z)  X_{\ell},
\qquad z \in \Bbb C^*,
$$
where $\psi (z)  $ in $\operatorname{SL}(V_{\ell})$ acts naturally on the space  
$\Bbb P^* (V_{\ell})$ of all hyperplanes in $V_{\ell}$ passing through the origin.
We then consider the map
$$
\pi : \mathcal{X}^{\psi} \to \Bbb A^1
$$ 
induced by the projection of $\Bbb A^1 \times \Bbb P^* (V_{\ell})$ to the first factor $\Bbb A^1$.
Moreover, for the hyperplane bundle $\mathcal{O}_{\Bbb P^*(V_{\ell})}(1)$ on $\Bbb P^*(V_{\ell})$, we put 
$$
\mathcal{L}^{\psi}\, :=\,\operatorname{pr}_2^*\mathcal{O}_{\Bbb P^*(V_{\ell})}(1)_{|\mathcal{X}^{\psi}}.
$$
Here $\operatorname{pr}_2 : \Bbb A^1 \times \Bbb P^* (V_{\ell}) \to \Bbb P^* (V_{\ell})$
denotes the projection to the second factor.
For the dual space $V_{\ell}^*$ of $V_{\ell}$,
the $\Bbb G_m$-action on $\Bbb A^1 \times V^*_{\ell}$ defined by
$$
\Bbb G_m \times (\Bbb A^1 \times V^*_{\ell})\to \Bbb A^1 \times V^*_{\ell},
\quad (t, (z, p))\mapsto  (tz, \psi (t) p),
$$
naturally induces a $\Bbb G_m$-action on $\Bbb A^1 \times \Bbb P^*(V_{\ell})$ and $\mathcal{O}_{\Bbb P^*(V_{\ell})}(-1)$, where $\operatorname{SL}(V_{\ell})$ acts on $V_{\ell}^*$ by the contragradient representation.
This induces $\Bbb G_m$-actions on $\mathcal{X}^{\psi}$ and $\mathcal{L}^{\psi}$, and $\pi : \mathcal{X}^{\psi} \to \Bbb A^1$ is a projective morphism with relative very ample line bundle 
$\mathcal{L}^{\psi}$ such that
$$
(\mathcal{X}^{\psi}_z, \mathcal{L}_z^{\psi})\; \cong \; (X,L^{\otimes \ell}),
\qquad z \neq 0,
$$
where $\mathcal{L}_z^{\psi}$ is the restriction of $\mathcal{L}^{\psi}$ to 
the fiber $\mathcal{X}^{\psi}_z := \pi^{-1}(z)$ over $z$.
Then a  triple $\mu = ({\mathcal{X}}, {\mathcal{L}}, \psi )$ is called 
a {\it test configuration for $(X,L)$}, if we have both $\mathcal{X} = \mathcal{X}^{\psi}$ and  $\mathcal{L}= \mathcal{L}^{\psi}$.
Note that $\ell$ is called the {\it exponent} of $({\mathcal{X}}, {\mathcal{L}}, \psi )$.
%Let
%$$
%\psi^{\operatorname{PGL}} : \Bbb G_m \to  \operatorname{PGL} (V_{\ell}).
%$$
%be the homomorphism induced by $\psi$. 
If $\psi$ %$\psi^{\operatorname{PGL}}$ 
is trivial, 
$({\mathcal{X}}, {\mathcal{L}}, \psi )$ is called {\it trivial}.
For the Lie algebra $\frak{sl}(V_{\ell})$, 
we now define a  symmetric 
bilinear form $\langle \;, \;\rangle^{}_{\ell}$
by
$$
\; \langle u, v\rangle^{}_{\ell} \; =\; \operatorname{Tr}(u v)/\ell^{n+2},
 \qquad u, v\in \frak{sl}(V_{\ell}),
$$
whose asymptotic limit as $\ell \to \infty$ plays an important role (cf.~\cite{Sz}) 
in the study of relative stability for test configurations.
Let $\frak t$ be the Lie algebra for $T$.
Since the infinitesimal action of $\frak t$
 on $X$ lifts to an infinitesimal bundle action of $\frak t$ on $L$, 
we can view $\frak t$ as a Lie subalgebra, denoted by $\frak t_{\ell}$, of $\frak{sl}(V_{\ell})$ by considering the traceless part.
Let
 $$
 u (\psi )\; :=\; \psi_*(t\partial/\partial t)\;\in\; \frak{sl}(V_{\ell})
 $$
be the generator for the infinitesimal action of $\operatorname{Lie}(\Bbb G_m )$ on $V_{\ell}$. 
We now consider the set $\mathcal{M}$ of all sequences $\{\mu_j \}$ of test configurations 
$$
\mu_j  \, =\, (\mathcal{X}_j, \mathcal{L}_j, \psi_j ), \qquad j =1,2,\dots,
\leqno{(2.2)}
$$
for $(X,L)$  such that $\langle u (\psi_j ), v\rangle^{}_{\ell} = 0$ for all $j$ and $v \in \frak t_{\ell}$,
and that the exponent $\ell_j$ of the polarized test configuration $\mu_j$ satisfies the following 
condition:
%\smallskip
%a) 
$$
\text{$\ell_j \to +\infty$,  \;\; as $j \to \infty$.}
$$
%\par
%b) $\ell_j | \ell_{j+1}$ for all positive integers $j$.
%\smallskip\noindent
In view of \cite{M}, to each $\{\mu_j\}\in  
\mathcal{M}$, we can assign the Donaldson-Futaki invariant 
$F_1 (\{\mu_j\}) \in  \Bbb R \cup \{-\infty\}$, see (5.8) below.
Then a strong version of relative
K-stability or relative K-semistability is defined as follows:

\medskip\noindent
{\em Definition\/ $2.3$}. 
(1) A polarized algebraic manifold $(X,L)$ is called 
{\it strongly $K$-semistable relative to $T$}, 
if $F_1(\{\mu_j\} ) \leq 0$ for all 
$\{\mu_j\} \in 
\mathcal{M}$.

\medskip\noindent
(2) A strongly K-semistable polarized algebraic manifold $(X,L)$
 is called 
{\it strongly $K$-stable relative to $T$},
if  for each $\{\mu_j\}\in\mathcal{M}$ satisfying $ F_1(\{\mu_j\} ) = 0$,
the test configurations $\mu_j$ are trivial for all sufficiently large $j$.

\medskip
Recall that these stability concepts are independent of the choice of the Hermitian 
metric $h$ for $L$ (see \cite{MN2}).

\section{Asymptotic relative Chow stability}

For the Lie subalgebra $\frak t_{\ell} \subset \frak{sl}(V_{\ell})$, let $\frak z_{\ell} := \{ u \in\frak{sl}(V_{\ell})\,;\, [u, \frak t_{\ell}] = 0\}$ be 
its centralizer in $\frak{sl}(V_{\ell})$.
Then by $\frak g_{\ell}$, we denote 
the set of all $u \in \frak z_{\ell}$ such that 
$$
\langle u, v\rangle^{}_{\ell}\; =\; 0, 
\qquad \text{ for all $ v\in \frak t_{\ell}$}.
$$
Let $T_{\ell}$ be the algebraic torus in $\operatorname{SL}(V_{\ell})$ with the Lie algebra $\frak t_{\ell}$. 
Then the centralizer $Z_{\ell}$ of $T_{\ell}$ in $\operatorname{SL}(V_{\ell})$ has the Lie algebra 
$\frak t_{\ell}$.
Let $(\frak g_{\ell})_{\Bbb Z}$ be the set of all $u$ in the kernel of the exponential map 
$$
\frak g_{\ell} \owns u \mapsto \exp (2\pi \sqrt{-1} u)\in Z_{\ell}
$$ 
such that the circle group $\{ \, \exp (2\pi s\sqrt{-1} u)\,;\, s \in \Bbb R\,\}$ acts isometrically 
on $(V_{\ell}, \rho_{\ell})$.
To each element $u$ in $(\frak g_{\ell})_{\Bbb Z}$, we assign a one-parameter group  
$\psi_u : \Bbb G_m \to \operatorname{SL}(V_{\ell})$, with the infinitesimal generator $u$, 
defined by
$$
\Bbb G_m \;\owns\; \exp (2\pi s\sqrt{-1}) \;\mapsto\;  \exp (2\pi s\sqrt{-1} u)\; \in\;
 \operatorname{SL}(V_{\ell}),
 \leqno{(3.1)}
$$
where $s$ runs through the set of all complex numbers.
Let $d_{\ell}$ denote the degree of $X_{\ell} := \Phi_{\ell}(X)$ in the projective space $\Bbb P^*(V_{\ell})$.
For each $u\in (\frak g_{\ell})_{\Bbb Z}$, let $\mathcal{X}^u$ denote the 
subvariety of $\Bbb A^1\times \Bbb P^*(V_{\ell})$ obtained as the closure of 
$$
\bigcup_{t\in \Bbb C^*}\; \{t\}\times \psi_{u}(t) X_{\ell}
$$
in $\Bbb A^1\times \Bbb P^*(V_{\ell})$, where $\operatorname{SL}(V_{\ell})$ 
acts naturally on the set $\Bbb P^*(V_{\ell})$ of all hyperplanes in $V_{\ell}$ passing through the origin.
Put $\mathcal{L}^u:= \operatorname{pr}_2^* \mathcal{O}_{\Bbb P^*(V_{\ell})}(1)$,
where $\operatorname{pr}_2: \mathcal{X}^u \to \Bbb P^*(V_{\ell})$ denotes the restriction to $\mathcal{X}^u$
of the projection to the second factor: $\Bbb A^1\times \Bbb P^*(V_{\ell})\to\Bbb P^*(V_{\ell})$.
Then for each $u\in (\frak g_{\ell})_{\Bbb Z}$, the triple
$$
\mu_u \; =\; (\mathcal{X}^u, \mathcal{L}^u, \psi_u ) 
\leqno{(3.2)}
%\qquad u\in (\frak t_{\ell}^{\perp})_{\Bbb Z},
$$
is called a {\it test configuration associated to $u$}.
For the fiber $\mathcal{X}^u_0$ of $\mathcal{X}^u$ over the origin in $\Bbb A^1$, 
we define $q(u)\in \Bbb Z$ by
$$
q(u): \;\text{the Chow weight for $\mathcal{X}^u_0$},
$$
where $\mathcal{X}^u_0$ is viewed as an algebraic cycle on $\{0\}\times \Bbb P^*(V_{\ell})
\, (\cong \Bbb P^*(V_{\ell}))$.
More precisely, $q(u)$ is  the weight at $\hat{\mathcal{X}}^u_0$ of the $\Bbb G_m$-action 
induced by $\psi_u$. Here
$$
\hat{\mathcal{X}}^u_0\;\in
\;W^*_{\ell}\,:= \;\{\operatorname{Sym}^{d_{\ell}}(V^*_{\ell})\}^{\otimes n+1}
$$
denotes the Chow form for $\mathcal{X}^u_0$, so that the associated element $[\hat{\mathcal{X}}^u_0]\in
\Bbb P^*(W_{\ell})$ is the Chow point for the cycle $\hat{\mathcal{X}}^u_0$ on $\Bbb P^*(V_{\ell})$.
Recall that (see \cite{MT1},\cite{MT2},\cite{MN1}):

\medskip\noindent
{\em Definition \/$3.3$.} (1) $(X,L^{\otimes\ell})$ is called {\it Chow-stable relative to $T$}, if 
the inequality $q_{\ell}(u) < 0$ holds for all 
$0 \neq u \in (\frak g_{\ell})_{\Bbb Z}$.

\smallskip\noindent
(2) $(X,L)$ is called {\it asymptotically Chow-stable relative to $T$}, if there exists a positive integer $\ell_0$ 
such that $(X,L^{\otimes \ell})$ is Chow-stable relative to $T$ for all positive integers $\ell$ 
satisfying $\ell \geq \ell_0$.

\medskip\noindent
{\em Remark \/$3.4$.}  In the preceding work \cite{MN1}, {\it relative Chow-stability} above
is referred to as {\it weak relative Chow-stability}.

\section{Polybalanced metrics and twisted Kodaira embeddings}

In this section, we consider the kernel $(\frak t_{\ell})_{\Bbb Z}$ of the exponential mapping 
$\frak t_{\ell} \owns y \mapsto \exp (2\pi \sqrt{-1}y) \in \operatorname{SL}(V_{\ell})$. 
Then by the infinitesimal $\frak t_{\ell}$-action on $V_{\ell}$, we can write the vector space $V_{\ell}$ as a direct sum of $\frak t_{\ell}$-eigenspaces
$$
V_{\ell}\; =\; \bigoplus_{k=1}^{m_{\ell}} \;V(\chi_{\ell;k} ),
$$
by setting $V(\chi_{\ell;k} ):= \{\sigma \in V_{\ell}\,;\,u\,\sigma = \chi_{\ell;k} (u) \sigma \text{ for all  $u\in (\frak t_{\ell})_{\Bbb Z}$}\}$ for mutually distinct  additive characters
$\chi_{\ell;k} \in \operatorname{Hom}( (\frak t_{\ell})_{\Bbb Z}, \Bbb Z)$, 
$ k=1,2,\dots,m_{\ell}$.  Put $n_{\ell, k} := \dim V(\chi_{\ell;k} )$.
Note that each $\chi_{\ell;k}$ naturally extends to an element of 
$\operatorname{Hom}( \frak t_{\ell}, \Bbb C)$, denoted also by $\chi_{\ell;k}$, by abuse of terminology.
Since $T_c$ acts isometrically on the K\"ahler manifold $(X,\omega )$, the Hermitian metric $h$ is also preserved by the $T_c$-action. Hence  
$$
V(\chi_{\ell;k} )\; \perp \; V_{\ell}(\chi_{\ell;k'} ), 
\qquad k\neq k',
$$
in terms of the Hermitian metric $\rho_{\ell}$ for $V_{\ell}$.
Consider the algebraic subgroup $S_{\ell}$ of $\operatorname{SL}(V_{\ell})$ defined by
$$
S_{\ell} \,:=\; \prod_{k=1}^{m_{\ell}}\;
\operatorname{SL}(V(\chi^{}_{\ell ;k} )),
$$
where the factor $\operatorname{SL}(V(\chi^{}_{\ell ;k} ))$ acts trivially
on $V(\chi^{}_{\ell ;k'} )$ if  $k \neq k'$.
Then the centralizer $Z(S_{\ell})$ of $S_{\ell}$ in $\operatorname{SL}(V_{\ell})$
consists of all diagonal matrices in $\operatorname{SL}(V_{\ell})$ acting on each $V(\chi^{}_{\ell ;k} )$ by constant scalar multiplication.
Hence 
$$
\frak z_{\ell}\; =\; \frak z (\frak s_{\ell}) \oplus \frak s_{\ell}
$$
where $\frak s_{\ell} := \operatorname{Lie}(S_{\ell})$ 
and $\frak z (\frak s_{\ell}) := \operatorname{Lie}(Z(S_{\ell}))$.
Let $\frak t_{\ell}^{\perp}$ be the orthogonal complement of $\frak t_{\ell}$ in $\frak z (\frak s_{\ell})$
defined by
$$
\frak t_{\ell}^{\perp}\; :=\;  \frak z (\frak s_{\ell}) \cap \frak g_{\ell}\; =\;\{\,u\in \frak z (\frak s_{\ell})\,;\,
\langle u, v\rangle^{}_{\ell} = 0\text{ for all $ v\in \frak t_{\ell}$}\,\},
$$
and the associated algebraic torus in $Z(S_{\ell})$ with the Lie algebra 
$\frak t_{\ell}^{\perp}$ will be denoted by $T_{\ell}^{\perp}$. 
Since $\frak g_{\ell} = \frak t_{\ell}^{\perp} + \frak s_{\ell}$, 
the associated reductive algebraic subgroup $G_{\ell}$ of $Z_{\ell}$ satisfies
$$
Z_{\ell} \; =\; T_{\ell} \cdot  G_{\ell}
\qquad\text{and}\qquad  G_{\ell} \; := \; T_{\ell}^{\perp} \cdot S_{\ell}.
\leqno{(4.1)}
$$
For integers $k$ and $\alpha$ with $1\leq k \leq m_{\ell}$ and $1 \leq \alpha \leq n_{\ell,k}$,
we put $i (k, \alpha ):= \alpha + \Sigma_{k'=1}^{k-1} n_{\ell,k'}$. 
Let $\mathcal{T} := \{\tau_1, \tau_2, \dots, \tau_{N_{\ell}}\}$ be a basis for the vector space $V_{\ell}$.

\bigskip\noindent
{\em Definition\/ $4.2$.}\,  (a) For each real number $r$, $O(\ell^{r})$ denotes a function $u$ or a $2$-form 
$\eta$ on $X$ satisfying $|u|\leq C_1 \ell^{r}$ or $ -C_1 \ell^{r} \omega \leq \eta \leq  C_1  \ell^{r}\omega$
for some positive real constant $C_1$ independent of the choice of $k$ and $\ell$.

\smallskip\noindent
(b) $\mathcal{T}$ is called {\it admissible\/} if for each $k$,
the space $V_{\ell} (\chi_k )$  admits a basis $\{\tau_{k,\alpha}\,;\, \alpha = 1,2, \dots, n_{\ell,k}\}$  such that

\smallskip\noindent
(b-1) \; $\tau_i \,\perp\, \tau_{i'}$\; by the Hermitian metric $\rho_{\ell}$, if $i \neq i'$;

\noindent
(b-2) \; $\tau_{i(k,\alpha )} \,=\, \tau_{k,\alpha}$, \quad $k=1,2,\dots,m_{\ell}; \; \alpha = 1,2,\dots, n_{\ell, k}$.

\smallskip\noindent
For an admissible basis $\mathcal{T}$ for $V_{\ell}$,
let $\mathcal{P}_{\mathcal{T}}$ denote the set of all elements in $Z_{\ell}$ written as positive real 
diagonal matrices in terms of the basis $\mathcal{T}$.
Then for every element $g$ of $\mathcal{P}_{\mathcal{T}}$, it is easily seen from (4.1) that
there exist $g' \in T_{\ell}\cap \mathcal{P}_{\mathcal{T}}$
 and $g''\in G_{\ell} \cap \mathcal{P}_{\mathcal{T}}$
such that
$\,g\,=\, g'\cdot g''$.

\bigskip\noindent
From now on until the end of this section, we assume that $(X, L^{\otimes \ell})$ is Chow-stable relative to $T$.
Then by \cite{MT2}, 
 the class $c_1(L)_{\Bbb R}$ admits a {\it polybalanced} metric $\omega_{\ell}$ in the following sense:
 There exist an admissible normal 
  basis $\mathcal{S} := \{\sigma_{i}\,;\,i =1,2, \dots, N_{\ell}\}$ for 
 $(V_{\ell}, \rho_{\ell})$ and an element  $v_{\ell}$ of $\frak t_{\ell}$
 such that the associated basis $\{\sigma_{k,\alpha}\,;\, \alpha = 1,2, \dots, n_{\ell,k}\}$ 
 for $V_{\ell} (\chi_k )$ satisfies
  
 \medskip\noindent
  (1) {\em  
  $\gamma_{\ell,k}:= 1 + \chi_k (v_{\ell})$ 
 is a positive real number such that $\omega_{\ell} = c_1(L;h_{\ell})$, where
 $h_{\ell}$ coincides with $(\Sigma_{k=1}^{m_{\ell}}\Sigma_{\alpha=1}^{n_{\ell, k}}
 \gamma_{\ell, k} |\sigma_{k,\alpha}|^2)^{-1/\ell}$ up to positive constant multiple}
 (cf.~\cite{MT2}, (3.7), (3.8), (3.9), (3.10), (3.11));

  \smallskip\noindent
 (2) 
 {\em For the Hermitian metric $V_{\ell} \owns \sigma', \sigma'' 
 \mapsto 
 \langle \sigma', \sigma''\rangle_{\hat{\rho}_{\ell}} := 
 \int_X (\sigma',\sigma'')_{h_{\ell}} \omega_{\ell}^n$ on $V_{\ell}$, 
the above $\mathcal{S}$ forms an orthonormal 
 basis for $(V_{\ell}, \hat{\rho}_{\ell})$} (cf.~\cite{MT2}, p.853);
 
 \smallskip\noindent
(3) {\em The $\ell$-th twisted Bergman kernel 
 $\hat{B}_{\ell}(\omega_{\ell}) :=  \Sigma_{k=1}^{m_{\ell}}\Sigma_{\alpha =1}^{n_{\ell, k}}\, \gamma_{\ell,k}|\sigma_{k,\alpha}|_{h_{\ell}}^2$ 
is the constant  function
$N_{\ell}/c_1(L)^n[X]$, where $N_{\ell}:= \dim V_{\ell}$} (cf.~\cite{MT2},  (3.9), (3.10)). 

%\medskip\noindent
%Then the $\ell$-th Bergman kernel $B_{\ell}(\omega_{\ell})
%:= \Sigma_{k=1}^{m_{\ell}}\Sigma_{\alpha =1}^{n_{\ell, k}}\, |\sigma_{k,\alpha}|_{h_{\ell}}^2$
% is written as 
%\begin{align*}
%&B_{\ell}(\omega_{\ell}) = \hat{B}_{\ell}(\omega_{\ell}) - 
%\Sigma_{k=1}^{m_{\ell}}\Sigma_{\alpha =1}^{n_{\ell, k}}\,
%( \gamma_{\ell,k}-1)\,|\sigma_{k,\alpha}|_{h_{\ell}}^2\\
%&= N_{\ell}/c_1(L)^n[X]\,\{\,1\,+\, 
%\hat{B}_{\ell}(\omega_{\ell})^{-1}\Sigma_{k=1}^{m_{\ell}}
%\Sigma_{\alpha =1}^{n_{\ell, k}}\,( \gamma_{\ell,k}-1)\,|\sigma_{k,\alpha}|_{h_{\ell}}^2\,\}
%\end{align*}

\smallskip\noindent
(4) {\em There exist uniformly bounded smooth functions $H_{\ell}$, $\ell \geq 1$, on $X$ such that
the $\ell$-th Bergman kernel $B_{\ell}(\omega_{\ell}):= \Sigma_{k=1}^{m_{\ell}}\Sigma_{\alpha =1}^{n_{\ell, k}}\, |\sigma_{k,\alpha}|_{h_{\ell}}^2$ is written as
 $$
 \{N_{\ell}/c_1(L)^n[X]\} + H_{\ell} \ell^{n-1} + O(\ell^{n-2}),
 $$
 where $H_{\ell}$ is a real Hamiltonian function on $X$,
so that
$(\sqrt{-1}/2\pi )\bar{\partial}H_{\ell} = i_{y_{\ell}}\omega_{\ell}$ for some holomorphic vector field $y_{\ell} = O(\ell^2 )v_{\ell}\in \frak t $ on $X$}
(cf~\cite{MT2}, (1.4)).

\bigskip\noindent
Let $\mathcal{T} := \{\tau_1, \tau_2, \dots, \tau_{N_{\ell}}\}$ be an admissible orthonormal basis for 
$(V_{\ell}, \rho_{\ell})$, and let $\{\tau^*_1, \tau^*_2, \dots, \tau^*_{N_{\ell}}\}$ denote its dual basis for 
the dual vector space $V_{\ell}^*$.
The associated Kodaira embedding $\Phi_{\ell} : X \hookrightarrow 
\Bbb P^{N_{\ell}-1}(\Bbb C ) $ is defined by
$$
\Phi_{\ell} (x) := (\tau_1 (x): \tau_2 (x): \dots : {\tau}_{N_{\ell}}(x)),
\qquad x\in X.
\leqno{(4.3)}
$$ 
The subset 
$\{\tau_{k,\alpha}\,;\, \alpha =1,2,\dots, n_{\ell,k}\}$  of $\mathcal{T}$ forms an orthonormal 
basis
for $(V_{\ell}(\chi_k ), \rho_{\ell})$, where we write each $\tau_{k,\alpha}$ as $\tau_{i(k, \alpha )}$ for simplicity.
We now put
$$
\hat{\tau}_{k,\alpha}\; :=\; \sqrt{\gamma_{\ell,k}}\, \tau_{k,\alpha}.
\leqno{(4.4)}
$$
Again by writing $\hat{\tau}_{k,\alpha}$ as $\hat{\tau}_{i(k,\alpha )}$, we have a basis 
$\{\hat{\tau}_i\,;\, i=1,2,\dots, N_{\ell}\}$ for $V_{\ell}$.
Then 
the embedding $\hat{\Phi}_{\ell} : X \hookrightarrow  
\Bbb P^{N_{\ell}-1}(\Bbb C ) =\{(z_1:z_2:\dots:z_{N_{\ell}})\}$ 
defined by
$$
\hat{\Phi}_{\ell} (x) := (\hat{\tau}_1 (x): \hat{\tau}_2(x):\, \dots \,: \hat{\tau}_{N_{\ell}}(x)),
\qquad x \in X,
\leqno{(4.5)}
$$
is called the {\it twisted Kodaira embedding\/} of $X$ for the complete linear system $|L^{\otimes \ell}|$
on $X$. Here $(z_1,z_2,\dots, z_{N_{\ell}})$ is the complex coordinates for the 
vector space $\Bbb C^{N_{\ell}}$ with standard 
Hermitian metric 
$$
\omega_{\operatorname{st}} \;:=\;
\sqrt{-1}\, ( dz_1\wedge  d\bar{z}_1 + d\bar{z}_2 \wedge d\bar{z}_2 + \dots + dz_{N_{\ell}}\wedge 
d\bar{z}_{N_{\ell}}),
$$
so that the complex projective space $\Bbb P^{N_{\ell}-1}(\Bbb C )$ above is identified with $\Bbb P^*(V_{\ell})$
via the isometry  $(\Bbb C^{N_{\ell}}, \omega_{\operatorname{st}}) \cong (V_{\ell}, \rho_{\ell})^*$ 
defined by
 $$
\Bbb C^{N_{\ell}} \owns\, (z_i,z_2,\dots,z_{N_{\ell}})\;\,\leftrightarrow\, \;
\Sigma_{i=1}^{N_{\ell}}z_i\tau^*_i\,\in V_{\ell}^*,
$$ 
and  $\omega_{\operatorname{FS}}:= (\sqrt{-1}/2\pi ) \partial\bar{\partial}\log (\Sigma_{i=1}^{N_{\ell}}|z_i |^2)$ is the corresponding Fubini-Study form on the projective space $\Bbb P^{N_{\ell}-1}(\Bbb C )$.

\medskip\noindent
{\em Remark\/ $4.6$.}\,  By \cite{MT2}, we have \,$\chi_k (v_{\ell}) = \gamma_{\ell,k} -1= O(\ell^{-1})$, see (4) above.

\medskip\noindent
{\em Remark\/ $4.7$.}  Let $\gamma_{\ell}$ be the element of $\operatorname{GL}( N_{\ell}, \Bbb C ) $ 
which sends each $ z  = (z_1,z_2,\dots, z_{N_{\ell}})\in \Bbb C^{N_{\ell}}$ to 
$\hat{z}  
= (\hat{z}_1, \hat{z}_2,\, \dots, \, \hat{z}_{N_{\ell}})
\in \Bbb C^{N_{\ell}}$, where $\hat{z}$ is defined by
$$
\hat{z}_{i (k,\alpha )}\; :=\; \sqrt{\gamma_{\ell, k}}\,
z_{i (k,\alpha )}, \;\;
\quad k =1,2,\dots,m_j;\, \alpha = 1,2,\dots, n_{\ell, k}.
$$ 
In view of $\hat{\Phi}_{\ell} = \gamma_{\ell}\circ \Phi_{\ell}$, we put
$\hat{\omega}_{\operatorname{FS}}:= (\sqrt{-1}/2\pi ) \partial\bar{\partial}
\log (\Sigma_{i=1}^{N_{\ell}}|\hat{z}_i |^2) = \gamma_{\ell}^*\omega_{\operatorname{FS}}$.
Then for ${\Phi}_{\ell}$ and $\hat{\Phi}_{\ell}$ in (4.3) and (4.5), we obtain
$$
\hat{\Phi}_{\ell}^*\omega_{\operatorname{FS}} \; =\; \Phi_{\ell}^*\hat{\omega}_{\operatorname{FS}}.
$$

\section{The twisted Donaldson-Futaki invariant}

Let $\mu_j = (\mathcal{X}_j, \mathcal{L}_j, \psi_j)$, $j = 1,2,\dots$,  be a sequence of test configurations 
for $(X,L)$ such that $\psi_j = \psi_{u_j}$ (cf.~(3.1)) 
for some $u_j \in (\frak g^{}_{\ell_j})^{}_{\Bbb Z}$, 
where
the exponent $\ell_j$ of $\mu_j$ satisfies 
$$
\ell_j \to \infty \quad \text{ as $j \to \infty$}.
$$
Put $d_j := \ell_j^{\,n} c_1(L)^n [X]$. Let $\mathcal{M}$ denote the set of all sequences $\{\mu_j \}$ of test configurations as above.
By abuse of therminology, $V_{\ell_j}$, $V_{\ell_j}(\chi_k )$, 
$N_{\ell_j}$, $\gamma_{\ell_j,k}$, $\gamma^{}_{\ell_j}$, $m_{\ell_j}$, $n_{\ell_j,k}$,
$\rho_{\ell_j}$, $\Phi_{\ell_j}$, $\hat{\Phi}_{\ell_j}$, $v_{\ell_j}$, $\frak t^{}_{\ell_j}$, 
$\frak g^{}_{\ell_j}$  will be written
simply as $V_j$, $V_j (\chi_k )$, $N_j$, $\gamma^{}_{j,k}$, $\gamma^{}_j$, $m_j$, $n_{j,k}$, $\rho_j$, $\Phi_j$, $\hat{\Phi}_j$, $v_j$, $\frak t^{}_j$, $\frak g_j$, respectively.
We then define
$$
W_{j} \; := \; \{\operatorname{Sym}^{d_j }(V_j )\}^{\otimes n+1},
$$
where $\operatorname{Sym}^{d_j}(V_j )$ denotes the $d_j$-th symmetric tensor product  of $V_j$.
The dual vector space $W_{j}^*$ of $W_{j}$ admits the Chow norm (cf.~\cite{Zh})
$$
W_{j}^* \owns w \;\mapsto \;\| w\|^{}_{\operatorname{CH}(\rho^{}_{j})} \in \Bbb R^{}_{\geq 0},
$$
associated to the Hermitian metric $\rho^{}_{j}$ on $V_j$ as in the introduction. Related to the 
twisted Kodaira embedding
$\hat{\Phi}_{j}: X \hookrightarrow  \Bbb P^{N_{\ell}-1}$ in (4.5) for the complete linear 
system $|L^{\otimes \ell_j}|$ on $X$, we consider the Chow form
$$
0 \neq \hat{X}_{j}\in W_{j}^*
$$ 
for the irreducible reduced algebraic cycle 
$\hat{\Phi}_{j}(X)$ on $\Bbb P^*(V_{j})$ such that the associated point $[\hat{X}_{j}]$ in $\Bbb P^*(W_j )$
is the Chow point for the cycle $\hat{\Phi}_{j}(X)$.
By setting $\hat{\mathcal{X}}_j := (\operatorname{id}_{\Bbb A^1}\times \gamma_j )_* (\mathcal{X}_j )$, we define a line bundle 
$\hat{\mathcal{L}}_j$ on $\hat{\mathcal{X}}_j$ by
$$
\hat{\mathcal{L}}_j \; := \; (\operatorname{id}_{\Bbb A^1}\times \gamma_j  )_* \mathcal{L}_j \;=\; \operatorname{pr}_2^*\mathcal{O}_{\Bbb P^*(V_j)}(1){}_{\,|\hat{\mathcal{X}}_j }.
$$
%For the homomorphism
% $\psi_{j}: T \to \operatorname{GL}(V_j)$,
%taking the real Lie subgroup 
%$$
%T_{\Bbb R}  = \{t \in T\,;\, t\in \Bbb R_+\}
%$$
%of the algebraic torus $T = \{t \in \Bbb C^*\}$, we define the associated lie group homomorphisms 
%$\psi^{\operatorname{SL}}_{j}: T \to \operatorname{SL}(V_j)$ by
%$$
%\psi^{\operatorname{SL}}_{j} (t) \; :=\; \frac{\psi_j (t)}{\;\det ( \psi_j (t))^{1/N_j}},
%\qquad t \in T_{\Bbb R}.
%$$
For the one-parameter group $\psi_j$, let $-\,b_{j;k, \alpha}$, 
$\alpha =1,2,\dots, n_{j,k}$,
be its weights
%, so that
%$\psi_j (t)\cdot \tau = t^{b_{j; k,\alpha }}\tau$ 
for the subspace
%all 
$%\tau \in 
V_j (\chi_k )$ of $V_j$, where $k = 1,2,\dots,m_{j}$.
We then have a basis $\{\,{\tau}_{k,\alpha}\,;\, \alpha = 1,2, \dots, n_{j,k}\,\}$
for $V_j (\chi_k )$ as in (4.4) such that 
$$
\psi_j (t)\cdot \hat{\tau}_{k,\alpha} \; =\; t_{}^{-\,b_{j;k, \alpha}}\, \hat{\tau}_{k,\alpha},
\qquad t \in \Bbb G_m.
$$
By identifying $\hat{\tau}_{k,\alpha}$ with $\hat{\tau}_{i(k,\alpha )}$, 
we consider the basis  $\{ \hat{\tau}^*_1, \hat{\tau}^*_2, \dots, \hat{\tau}^*_{N_j}\}$ for $V_j^*$
dual to $\{\hat{\tau}_1, \hat{\tau}_2, \dots, \hat{\tau}_{N_j}\}$. By writing $b_{j;k, \alpha}$ 
as $b_{j;i (k, \alpha )}$, we obtain
$$
\psi_j (t)\cdot \hat{\tau}^*_{i} \; =\; t_{}^{b_{j;i}}\,\hat{\tau}^*_{i},
\qquad i = 1,2, \dots, N_j.
$$ 
Each ${z} = ({z}_1, {z}_2, \dots , {z}_{N_j}) \in \Bbb C^{N_j}\setminus \{0\} $ sitting over 
$ [{z}]=({z}_1:{z}_2: \dots :{z}_{N_j}) \in \Bbb P^{N_j-1}(\Bbb C ) = \Bbb P^*(V_j)$ is written as 
$\Sigma_{i=1}^{N_j} {z}_i \hat{\tau}^*_{i}$, so that the action by $t \in \Bbb G_m$ 
on $\Bbb C^{N_j}$ is expressible in the form
$$
{z}= ({z}_1, {z}_2, \dots, {z}_{N_j}) \; \mapsto \; 
\psi_j (t)\cdot {z} =  (t^{b_{j;1}}{z}_1, t^{b_{j;2}}{z}_2, \dots, 
t^{b_{j;N_j}}{z}_{N_j}).
$$
For each $p\in\hat{\mathcal{X}_j }{}_{|t = 1} \,(= \hat{\Phi}_j (X))$,
writing $p$ as $\hat{\Phi}_j (x)$ for some $x\in X$,  we have a line $l_p$ through the origin in $\Bbb C^{N_j}$
associated to  $p = \hat{\Phi}_j (x)$ in $\Bbb P^{N_j-1}(\Bbb C )$. 
Let $t \in \Bbb G_m$. Put $p' := \psi_j (t)\cdot p \in \hat{\mathcal{X}}_t$, 
where we view $\hat{\mathcal{X}}_t$ in $ \{t\}\times \Bbb P^{N_j}(\Bbb C)$  
as a subvariety of $\Bbb P^{N_j}(\Bbb C)$ by identifying $ \{t\}\times \Bbb P^{N_j}(\Bbb C)$ 
with $ \Bbb P^{N_j}(\Bbb C)$.
Let $l_{p'}$ denote the line through the origin in $\Bbb C^{N_j}$ associated to
$p'$ in $\Bbb P^{N_j}(\Bbb C)$.
Then 
$$
l_p \to l_{p'}, \qquad {z} \mapsto \psi_j (t) \cdot {z},
$$
naturally defines a map of $\hat{\mathcal{L}}_p^{-1}$ onto $\hat{\mathcal{L}}_{p'}^{-1}$.
Since $\psi_j (t)$, $t \in \Bbb G_m$, and $\gamma_j$ (cf.~Remark 4.7) commute,
the $\Bbb G_m$-action on ${\mathcal{X}}_j$ naturally induces a 
$\Bbb G_m$-action on $\hat{\mathcal{X}}_j$, while the $\Bbb G_m$-action on $\hat{\mathcal{X}}_j$
lifts to a $\Bbb G_m$-action on $\hat{\mathcal{L}}^{-1}$ as above.  
 %let $l_p$ denote the fiber of 
%$L^{-\ell_j}\, (= \mathcal{L}_j^{-1}{}_{|\hat{\Phi}_j (X)$ over $p$.
%Then the set of all $a = (a_1,a_2, \dots ,a_{N_j}) \in \Bbb C^{N_j}\setminus \{0\}$ such that 
 %$[a] = p$ is identified with $l_p \setminus \{0\}$.
Put $\|b_j\|^{}_1:= \Sigma_{k=1}^{m_j} \Sigma_{\alpha =1}^{n_{j,k}}\,  |b^{}_{j;k,\alpha}|$.
Then we define $\|\psi_j\|_1$ and $\|\psi_j\|_{\infty}$  by
$$
\begin{cases}
&\|\psi_j\|^{}_1 \,:= \;\|b_j\|^{}_1/\ell_j^{n+1}\,=\;
\Sigma_{k=1}^{m_j} \Sigma_{\alpha =1}^{n_{j,k}}\,  |b^{}_{j;k,\alpha}|/\ell_j^{n+1},\\
%\quad\; \text{ and }\; \quad 
&\|\psi_j\|^{}_{\infty} :=\; \max \{\, |b^{}_{j;k,\alpha}|/\ell_j \,;\,
1 \leq  k\leq m_j,\,1\leq  \alpha \leq  n_{j,k}\,\}.
\end{cases}
%\leqno{(5.1)}
$$
Let $\delta (\psi_j )$ denote $\|\psi_j\|_{\infty}/\|\psi_j\|_1$ or $1$ according as 
$\|\psi_j\|_{\infty} \neq 0$ or $\|\psi_j\|_{\infty}=0$. 
If $\|\psi_j\|_{\infty} \neq 0$, we write $t\in \Bbb R_+$ 
as 
$$
t = \exp (s/\|\psi_j\|_{\infty}),
\leqno{(5.1)}
$$
while we require no relations between $s\in \Bbb R$ and $t\in T_{\Bbb R}$
if $\|\psi_j\|_{\infty}=0$. 
Since $\operatorname{SL}(V_j)$ acts naturally on  $W_j^*$, 
%by writing $t\in T_{\Bbb R}$ 
%as $t = \exp (s/\|\psi_j\|_{\infty})$,
we define a function $\hat{f}_j\, =\, \hat{f}_j (s)$ on $\Bbb R$ by 
$$
\hat{f}_j (s)\; :=\; \delta (\psi_j )\,\ell_j^{-n}\log 
\| \psi_j (t)\cdot \hat{X}_j \|^{}_{\operatorname{CH}(\rho^{}_{j})},
\qquad s \in \Bbb R.
%\leqno{(5.3)}
$$
To understand the action of $\Bbb G_m$ on $(\mathcal{X}, \mathcal{L}_j)$ clearly, 
we consider 
We shall now show the following:

\medskip\noindent
{\bf Proposition 5.2.} {\em If $j \gg 1$, then the derivative 
$(d\hat{f}_j /ds)_{|s=0}$ is bounded from above  by a positive real constant $C_2$ independent of the 
choice of $j$.}

\medskip\noindent
{\em Proof}:
If $\|\psi_j\|_{\infty} = 0$, then $\hat{f}_j$ is a constant function, so that $(d\hat{f}_j /ds)_{|s=0} = 0$. 
Hence we may assume $\|\psi_j\|_{\infty} \neq  0$ without loss of generality. For the K\"ahler metric $\omega = c_1(L;h)$ on $X$, 
let $\hat{B}_j (\omega )$ denote the $\ell_j$-th asymptotic twisted Bergman kernel defined by 
$$
\hat{B}_j (\omega )\; :=\; \Sigma_{i =1}^{N_j}\; |\hat{\tau}_i |^2_h,
\leqno{(5.3)}
$$
where we put $|\hat{\tau}_i |^2_h := (\hat{\tau}_i, \hat{\tau}_i )_h$ for the basis
$\{\hat{\tau}_i\,;\, i=1,2,\dots, N_j\}$ for $V_j$ as in (4.4) and (4.5).
For the Fubini-Study form $\omega_{\operatorname{FS}}:= (\sqrt{-1}/2\pi )\partial\bar{\partial}\log 
\Sigma_{i =1}^{N_j}|z_i|^2$ on $\Bbb P^{N_j -1}(\Bbb C )$, by \cite{Ze}, we see from 
(4.4) and Remark 4.6 that
\begin{align*}
& |\hat{B}_j - (1/n! )\ell_j^n | \; \leq \; C_3 \ell_j^{n-1},
\tag{5.4}\\
& (\ell_j - C_4)\omega \;\leq \; \hat{\Phi}^*_j \omega_{\operatorname{FS}}
\;\leq \; (\ell_j + C_4 )\omega,
\tag{5.5}
\end{align*}
for some positive real constants $C_3$ and $C_4$ independent of $j$,
where for (5.5), see for instance  \cite{MT}, p.234.
By \cite{Zh} (see also \cite{MT1}), $(d\hat{f}_j /ds)_{|s=0}\,$ is 
$$
(\|\psi_j \|^{}_1 \ell_j^n )^{-1}\int_X \hat{B}_j(\omega )^{-1}
(\Sigma_{k=1}^{m_j}\Sigma_{\alpha =1}^{n_{j,k}}\, b_{j;k,\alpha} |\hat{\tau}_{k,\alpha} |_h^2)\, \hat{\Phi}^*_j\omega^n_{\operatorname{FS}}.
\leqno{(5.6)}
$$
Since $u_j$ belongs to $(\frak g_{j})_{\Bbb Z}$, and since $v_j$ belongs to $\frak t_{j}$, 
in view of the equality %$\chi_k (u_j) = b_{j,k}$ and 
$\chi_k (v_j ) = \gamma_{j,k} -1$,
we obtain
$$
0\, =\,\langle u_j, v_j \rangle_{\ell_j } \, = \, -\Sigma_{k=1}^{m_j}\Sigma_{\alpha =1}^{n_{j,k}}\, b_{j;k,\alpha} (\gamma_{j,k} -1) 
\, =\, -\Sigma_{k=1}^{m_j}\Sigma_{\alpha =1}^{n_{j,k}}\, b_{j;k,\alpha}\gamma_{j,k},
$$
so that $\int_X 
\Sigma_{k=1}^{m_j}\Sigma_{\alpha =1}^{n_{j,k}} b_{j;k,\alpha} |\hat{\tau}_{k,\alpha} |_h^2 \omega^n =
\Sigma_{k=1}^{m_j}\Sigma_{\alpha =1}^{n_{j,k}}b_{j;k,\alpha}\gamma_{j,k} = 0$. By setting
$$
\begin{cases}
& R_1:= \int_X  \hat{B}_j (\omega )^{-1}(\Sigma_{k =1}^{m_j} \Sigma_{\alpha =1}^{n_{j,k}}
b_{j;k,\alpha}|\hat{\tau}_{k,\alpha} |_h^2)\{\hat{\Phi}^*_j\omega^n_{\operatorname{FS}} - (\ell_j \omega )^n\},
\\
&R_2 := \int_X \{ \hat{B}_j (\omega )^{-1} - (\ell_j /n! )^{-1}\}(\Sigma_{k =1}^{m_j} \Sigma_{\alpha =1}^{n_{j,k}}
b_{j;k,\alpha}|\hat{\tau}_{k,\alpha} |_h^2)(\ell_j \omega )^n,
\end{cases}
$$
we see from (5.6) that $(d\hat{f}_j /ds)_{|s=0} = (\|\psi_j \|^{}_1 \ell_j^n )^{-1}(R_1 + R_2)$,
while by (5.4) and (5.5), if $j \gg 1$, $(\|\psi_j \|^{}_1 \ell_j^n )^{-1}R_1 $ is less than or equal to
%(\|\psi_j \|^{}_1 \ell_j^n )^{-1}R_1 \leq 
\begin{align*}
&\ell_j \|b_j\|_1^{-1}
\int_X  \hat{B}_j (\omega )^{-1}
(\Sigma_{k =1}^{m_j} \Sigma_{\alpha =1}^{n_{j,k}}
|b_{j;k,\alpha}|\,|\hat{\tau}_{k,\alpha} |_h^2)\, |\hat{\Phi}^*_j\omega^n_{\operatorname{FS}} - (\ell_j \omega )^n|\\
&\leq \;  C_4 \ell_j \|b_j\|_1^{-1}\int_X (\ell^n_j /n! )^{-1}
(\Sigma_{k =1}^{m_j} \Sigma_{\alpha =1}^{n_{j,k}}
|b_{j;k,\alpha}|\,|\hat{\tau}_{k,\alpha} |_h^2)\,\ell_j^{n-1}\omega^n \\
&= \; n! \,C_4(\Sigma_{k=1}^{m_j}\Sigma_{\alpha =1}^{n_{j,k}} |b_{j;k,\alpha}|)^{-1}\int_X (\Sigma_{k =1}^{m_j} \Sigma_{\alpha =1}^{n_{j,k}}
|b_{j;k,\alpha}|\,|\hat{\tau}_{k,\alpha} |_h^2)\,\omega^n \; \leq  \; 2\,n!\,C_4,
\end{align*}
for a positive real constant $C_4$ independent of $j$, where the last inequality follows from
(4.4) and Remark 4.6.
 Finally, if $j \gg 1$, by (5.4)
\begin{align*}
&(\|\psi_j \|^{}_1 \ell_j^n )^{-1} R_2  \leq  C_5\ell_j \|b_j\|_1^{-1}
\int_X \ell_j^{-n-1}(\Sigma_{k =1}^{m_j} \Sigma_{\alpha =1}^{n_{j,k}}
b_{j;k,\alpha}|\hat{\tau}_{k,\alpha} |_h^2)(\ell_j \omega )^n\\
&\leq C_5(\Sigma_{k=1}^{m_j}\Sigma_{\alpha =1}^{n_{j,k}} |b_{j;k,\alpha}|)^{-1}\int_X (\Sigma_{k =1}^{m_j} \Sigma_{\alpha =1}^{n_{j,k}}
|b_{j;k,\alpha}|\,|\hat{\tau}_{k,\alpha} |_h^2)\,\omega^n \leq  2 C_5,
\end{align*}
for a positive real constant $C_5$ independent of $j$, where in the last inequality we again used (4.4) 
and Remark 4.6. Thus we obtain $(d\hat{f}_j /ds)_{|s=0} = (\|\psi_j \|^{}_1 \ell_j^n )^{-1}(R_1 + R_2)
\leq 2\,n!\,C_4 + 2\, C_5$, as required.

\medskip
By the boundedness $(d\hat{f}_j /ds)_{|s=0} \leq C$ and
by the convexity of the function $\hat{f}_j = \hat{f}_j (s)$, we easily see that $\varliminf_{j \to \infty} d\hat{f}_j /ds$ is a well-defined non-decreasing function
in $s$. 
We can now define the {\it twisted Donaldson-Futaki invariant} 
$\hat{F}_1 (\{\mu_j\}) \in \Bbb R \cup \{-\infty\}$ for $\{\mu_j \}\in\mathcal{M}$ by
$$
\hat{F}_1 (\{\mu_j\}) := \lim_{s\to -\infty}\{\varliminf_{j\to \infty} d\hat{f}_j/ds \}.
\leqno{(5.7)}
$$
Next by setting $\ell = \ell_j$ in (4.3), 
we consider the Kodaira embedding $\Phi_j$, and view
the image $\Phi_j (X)$ as an irreducible reduced
algebraic cycle on the projective space $\Bbb P^*(V_j)$. 
We then consider its Chow form
$$
0 \neq X_j \in W_j^*,
$$
so that the corresponding point $[X_j]$ in $\Bbb P^*(W_j )$ is the Chow point for the cycle $\Phi_j (X)$.
We now define a function $f_j = f_j (s)$ on $\Bbb R$ by 
$$
f_j (s) \;:=\; \delta (\psi_j ) \ell_j^{-n}\log \|\psi_j (t)\cdot X_j \|_{\operatorname{CH}(\rho_j )},
\qquad s \in \Bbb R,
$$
where $s\in \Bbb R$ and $t \in \Bbb R_+$ satisfy (5.1) 
or have no relation according as $\|\psi_j \|_{\infty} \neq 0$ or $\|\psi_j \|_{\infty} = 0$.
Recall that,  in \cite{M},  we defined the {\it Donaldson-Futaki invariant} $\hat{F}_1 (\{\mu_j\}) \in \Bbb R \cup \{-\infty\}$ for $\{\mu_j \}\in\mathcal{M}$ by
$$
{F}_1 (\{\mu_j\}) := \lim_{s\to -\infty}\{\varliminf_{j\to \infty} d{f}_j/ds \}.
\leqno{(5.8)}
$$
Put $\Delta := \varlimsup_{j\to \infty} \, \{\delta (\psi_j )/\ell_j\}$.
We now claim the following:

\medskip\noindent
{\bf Lemma 5.9.}\qquad
{\em $\hat{F}_1 (\{\mu_j\}) =  F_1 (\{\mu_j\})$ for all $\{\mu_j \} \in \mathcal{M}$.}

\medskip\noindent
{\em Proof\/}: For the case $\Delta > 0$, see \cite{M2} and \cite{MN2}.
Hence we assume that $\Delta = 0$.
%As in Remark 4.7, let $\gamma_j\in\operatorname{GL}(V_j)$ be such that
% $$
% \gamma_j \tau_i  := \hat{\tau}_i, 
% \qquad i=1,2,\dots,N_j,
% $$
%by the notation in (4.3) and (4.4) applied to $\ell = \ell_j$.
If $\|\psi_j \|_{\infty} =0$, then on the whole real line $\Bbb R$ 
$df_j/ds = \; 0 \;=\;  d\hat{f}_j/ds$,
and hence we may assume that 
$$
\|\psi_j \|_{\infty}\; \neq \; 0.
$$
For simplicity, we put $\ln (\theta ) := (1/2\pi )\log \theta$ for Hermitian metrics $\theta$ 
for $L^{-1}$.
Put $\theta_{\operatorname{FS}} := \Sigma_{i=1}^{N_j} |z_i |^2$ and 
$\hat{\theta}_{\operatorname{FS}} = \Sigma_{i=1}^{N_j} |\hat{z}_i |^2$,
where $\hat{z}_i$ is as in Remark 4.7.
We now consider the associated Fubini-Study forms on $X$ defined by
$$
\omega_{\operatorname{FS}} :=  \sqrt{-1}\,{\partial}\bar{\partial}\ln \theta_{\operatorname{FS}}{}^{}_{\,|z\, =\, \tau }\quad \text{ and }\quad
\hat{\omega}_{\operatorname{FS}} := \sqrt{-1}\,
 {\partial}\bar{\partial}\ln \hat{\theta}_{\operatorname{FS}}{}^{}_{\,|\hat{z}\, =\, \hat{\tau} }
$$
 where the equality $z\,=\,\tau$ means that 
$z_i = \tau_i $, for all $i=1,2,\dots,N_j$, and the equality $\hat{z}\,=\,\hat{\tau}$ means that 
$\hat{z}_i = \hat{\tau}_i $, for all $i=1,2,\dots,N_j$.
We then define $\theta_j(s) := \{\psi_j (t) \cdot \theta_{\operatorname{FS}}\}_{|z=\tau}$ and
$\hat{\theta}_j(s) := \{\psi_j (t) \cdot \hat{\theta}_{\operatorname{FS}}\}_{|\hat{z}=\hat{\tau}}$ by
\begin{align*}
& \theta_j(s) :=  (\Sigma_{i=1}^{N_j} |\psi_j (t)\cdot z_i |^2){\,}_{|z=\tau}
= \Sigma_{k=1}^{m_j}\Sigma_{\alpha =1}^{n_{j,k}}|t|^{2b_{j;k,\alpha}}|\tau_{k,\alpha}|^2,
\tag{5.10}
\\
&\hat{\theta}_j(s):= (\Sigma_{i=1}^{N_j} |\psi_j (t)\cdot\hat{z}_i |^2){\,}_{|\hat{z}=\hat{\tau}}
= \Sigma_{k=1}^{m_j}\Sigma_{\alpha =1}^{n_{j,k}} \gamma_{j,k}|t|^{2b_{j;k,\alpha}}|\tau_{k,\alpha}|^2.
\tag{5.11}
\end{align*}
Let $\Theta$ be the set of all Hermitian metrics $\theta$ for $L^{-1}$ such that 
the associated forms $\sqrt{-1}\,\partial\bar{\partial}\ln ( \theta )$ are K\"ahler.
Recall that, given elements $\theta'$, $\theta''$ in $\Theta$, an energy
$\kappa : \Theta\times \Theta \to \Bbb R$ is defined by
$$
\kappa (\theta', \theta'')\;:=\; \int^{s\,=\,a''}_{s\,=\,a'}\left ( \int_X  \{\partial \ln (\theta_s )/\partial s\}\,
\{ \sqrt{-1}\,\partial\bar{\partial}\ln ( \theta_s ) \}^n\right ) ds
\leqno{(5.12)}
$$
for a smooth path $\theta_s$, $a'\leq s \leq a''$, in $\Theta$ satisfying $\theta_{a'} = \theta'$ and $\theta_{a''}=\theta''$, where 
$\kappa (\theta', \theta'')$ is independent of the choice of the path. Then by \cite{Zh} (see also \cite{MT1} and \cite{Sn}), for each positive integer $p$,
we obtain
\begin{align*}
&f_j(1-p)  - f_j (-p)\; =  \;\delta (\psi_j )\,\ell_j^{-n}\,\kappa (\theta_j(-p), \theta_j (1-p)),
\tag{5.13}
\\
&\hat{f}_j(1-p)  - \hat{f}_j (-p)\; =  \;\delta (\psi_j )\,\ell_j^{-n}\,
\kappa (\hat{\theta}_j(-p), \hat{\theta}_j (1-p)).
\tag{5.14}
\end{align*}
By the mean value theorem, there exist real numbers $s_{j,p}$ and $\hat{s}_{j,p}$ on the open interval
$(-p, 1-p)$ such that
$$
\begin{cases}
&f_j(1-p)  - f_j (-p) = (df_j/ds)_{|s\, =\, s_{j,p}},\\
& \hat{f}_j(1-p)  - \hat{f}_j (-p) = (d\hat{f}_j/ds)_{|s\, =\, \hat{s}_{j,p}}. 
\end{cases}
\leqno{(5.15)}
$$
By subtracting (5.14) from (5.13), we see from (5.15) that
$$
\begin{cases}
&(df_j/ds)_{|s\, =\, s_{j,p}} \,-\, (d\hat{f}_j/ds)_{|s\, =\, \hat{s}_{j,p}} \\
&= \; \delta (\psi_j )\,\ell_j^{-n}\,\{\,\kappa (\theta_j(-p), \theta_j (1-p)) \,-\,
\kappa (\hat{\theta}_j(-p), \hat{\theta}_j (1-p))\,\}\\
&= \; \delta (\psi_j )\,\ell_j^{-n}\,\{\,\kappa (\theta_j(-p), \hat{\theta}_j (-p)) \,-\,
\kappa ({\theta}_j (1-p), \hat{\theta}_j(1-p))\,\}.
\end{cases}
\leqno{(5.16)}
$$
In view of (5.12), by (5.10) and (5.11), it follows from Remark 4.6 that
$$
|\,\kappa (\theta_j (s), \hat{\theta}_j (s))\,| \, \leq \, C_6\,\ell_j^{n-1}
\leqno{(5.17)}
$$ 
for some positive real constant $C_6$ independent of $s$ and $j$. Given an arbitrary positive real number 
$\varepsilon$,
since both $df_j/ds$ and $d\hat{f}_j/ds$ are non-decreasing functions on $\Bbb R$, (5.16) and (5.17) imply
\begin{align*}
&(df_j/ds)_{|s\, =\,1-p} \,-\, (d\hat{f}_j/ds)_{|s\, =\, -p}\;\geq \;
(df_j/ds)_{|s\, =\, s_{j,p}} \,-\, (d\hat{f}_j/ds)_{|s\, =\, \hat{s}_{j,p}}\\
&\geq \;-\,2\,\delta (\psi_j ) \,C_6 \, \ell^{-1}_j \;\geq \; -\,2\, C_6 \, (\Delta + \varepsilon )
\; = \; -\,2\,C_6\, \varepsilon,
 \end{align*}
for all $j \gg 1$ and all positive integers $p$. 
Now, let $j \to \infty$. Then
$$
\varliminf_{j \to \infty}  (df_j/ds)_{|s\, =\, 1-p}  \, - \, 
\varliminf_{j \to \infty}  (d\hat{f}_j/ds )_{|s \,=\,-p} \;\geq\;  -\,2\,C_6\, \varepsilon.
$$
Since both $df_j/ds$ and $d\hat{f}_j/ds$ are non-decreasing, by letting $p \to \infty$, we have
the inequality $F_1 (\{\mu_j \}) - \hat{F} (\{\mu_j \}) \; \geq\;  -\,2\,C_6\,  \varepsilon$.
Since $\varepsilon >0$ is arbitrary, further by letting $\varepsilon \to 0$, we obtain
$$
F_1 (\{\mu_j \}) \; \geq \;  \hat{F} (\{\mu_j \}).
\leqno{(5.18)}
$$
Next by using (5.16) and (5.17), we similarly 
have the following for every positive real number $\varepsilon $:
$$
(d\hat{f}_j/ds)_{|s\, =\,1-p} \,-\, (d{f}_j/ds)_{|s\, =\, -p}\; \geq \; -\,2\,C_6\,  \varepsilon,
$$
for all $j \gg 1$ and all positive integers $p$. Hence by letting $j \to \infty$, we obtain
$$
\varliminf_{j \to \infty}  (d\hat{f}_j/ds)_{|s\, =\, 1-p}  \, - \, 
\varliminf_{j \to \infty}  (d{f}_j/ds )_{|s \,=\,-p} \;\geq\;  -\,2\,C_6\,  \varepsilon.
$$
Let $p \to \infty$ and $\varepsilon \to 0$. Then we see that
$\hat{F}(\{\mu_j \}) \geq  F_1 (\{\mu_j \})$.
From this and (5.18), we now conclude that 
$\hat{F}(\{\mu_j \}) = F_1 (\{\mu_j \})$,
as required.

\section{Test configurations associated to polybalanced metrics}

Hereafter, we assume that the polarized algebraic manifold $(X,L)$ is strongly K-stable relative to $T$.
%while we fix  a Hermitian metric $h\in\mathcal{H}$ so that
%$\omega := c_1 (L;h)$ is K\"ahler.
Then by \cite{MN1},
$(X,L)$  is asymptotically Chow-stable relative to $T$, so that %for each $\ell\gg 1$, i.e., 
there exists an integer $\ell_0 \gg 1$ such that
for all $\ell $ with $\ell \geq \ell_0$, $(X,L)$  is Chow-stable relative to $T$.
Throughout this section, we assume  $\ell \geq \ell_0$, so that by \cite{MT2}, there exists
a polybalanced K\"ahler metric $\omega_{\ell} := c_1 (L; h_{\ell})$ in the sense that, 
by using the notation in Section 4, we have
$$
\Sigma_{i=1}^{N_{\ell}} \, |\hat{\sigma}_{i}|^2_{h_{\ell}} \; =\; N_{\ell}/c_1 (L)^n [X], 
\leqno{(6.1)}
$$
where $\hat{\sigma}^{}_{i (k,\alpha )} = \hat{\sigma}^{}_{k,\alpha } 
= \sqrt{\gamma_{\ell,k}}\,\sigma^{}_{ k,\alpha } = \sqrt{\gamma_{\ell,k}}\,\sigma^{}_{ i (k,\alpha )}$ 
as in (4.4)
for an admissible orthonormal basis 
$\{\sigma_1, \sigma_2, \dots,\sigma_{N_{\ell}}\}$ for 
$(V_{\ell}, \hat{\rho}_{\ell})$.
By comparing the Hermitian vector spaces 
$(V^{}_{\ell},\hat{\rho}^{}_{\ell})$ with $(V^{}_{\ell},\rho^{}_{\ell})$, we have 
admissible orthonormal bases
$$
\{\sigma^{}_{1}, \sigma^{}_{2}, \dots, \sigma^{}_{N_{\ell}}\}
\;\;\text{ and }\;\; \{\tau^{}_{1}, \tau^{}_{2}, \dots, \tau^{}_{N_{\ell}}\}
$$
for $(V^{}_{\ell},\hat{\rho}^{}_{\ell})$ and $(V^{}_{\ell},\rho^{}_{\ell})$, respectively, such that
there exit positive real numbers $ \lambda^{}_{\ell, i}$, $i =1,2,\dots,N_{\ell}$, satisfying
$$
\sigma^{}_{i} \; =\; \lambda^{}_{\ell, i}\,\tau^{}_{i},
\qquad i = 1,2,\dots,N_{\ell}.
\leqno{(6.2)}
$$
Replacing $h_{\ell}$ by a suitable positive constant, if necessary, 
we may assume that 
$\Pi_{i =1}^{N_{\ell}}\, \lambda^{}_{\ell, i} = 1$.
Define an element $\lambda_{\ell}$ of $Z_{\ell}$ by
$$
\lambda_{\ell} (\tau_i ) \; =\; \lambda^{}_{\ell, i}\, \tau_i,
\qquad i = 1,2,\dots, N_{\ell}.
$$
In view of (4.1), we have $Z_{\ell} = T_{\ell}\cdot G_{\ell}$, and hence 
for some $g \in T_{\ell}$ and $\lambda'_{\ell} \in G_{\ell}$
we can write $\lambda_{\ell}$ in the form
$$
\lambda_{\ell}\; = \;g\cdot \lambda'_{\ell} \;=\; \lambda'_{\ell}\cdot g,
$$
where both $g$ and $\lambda_{\ell}'$ are written as positive real diagonal matrices in terms of the basis 
$\{\tau_1,\tau_2,\dots,\tau_{N_{\ell}}\}$ (see the remark in Definition 4.2).
In view of (6.1) above, $g^*\omega_{\ell}$ is also polybalanced with the same 
weight $\gamma_{\ell, k}$ (cf.~(1) of Section 4), where we observe that the $T_{\ell}$-action
on $V_{\ell}$ is induced by the $T$-action on $X$. 
Put $g\cdot h_{\ell} := (g^{-1})^*h_{\ell}$ and $g\cdot \omega_{\ell}:= (g^{-1})^*\omega_{\ell}$.
Replacing the metrics $h_{\ell}$ and $\omega_{\ell}$
by the metrics $g\cdot h_{\ell}$ and $g\cdot \omega_{\ell}$, respectively, we may assume that $\lambda_{\ell} = \lambda'_{\ell}$,
and therefore
$$
\lambda_{\ell} \,\in\, G_{\ell},
\leqno{(6.3)}
$$
where (4) in Section 4 is true even if $(h_{\ell}, \omega_{\ell} )$ is replaced by 
$(g\cdot h_{\ell}, g\cdot \omega_{\ell})$.
Let $\{\underline{v}_{r}\,;\, r= 1,2,\dots, r_0 \}$ be the free basis for 
the free $\Bbb Z$-module $(\frak t_{\ell})_{\Bbb Z}$. 
Since $\tau_{k,\alpha}= \tau_{i(k,\alpha )}$
as in  (b-2) of Definition 4.2, by identifying $\Sigma_k\Sigma_{\alpha}\,  
\Bbb R \tau_{k,\alpha} \otimes \tau_{k,\alpha}^*$
with $\Bbb R^{N_{\ell}}$, we see that the intersection of $\frak g_{\ell}$ with this 
space $\Bbb R^{N_{\ell}}$ is 
$$
\frak S_{\ell} :=\{ \Sigma_{k}\Sigma_{\alpha}
\,  b_{k,\alpha }\tau_{k,\alpha} \otimes \tau_{k,\alpha}^*\in  
\Bbb R^{N_{\ell}}\, ;\,\Sigma_{k}\Sigma_{\alpha}
b_{k,\alpha }\chi^{}_{\ell;k} (\underline{v}_{r} ) = 0 \text{ for all $r$}\},
$$
where in the summation $\Sigma_k\Sigma_{\alpha}$, $k$ runs through $\{1,2,\dots,m_j\}$
and $\alpha$ runs through $\{1,2,\dots n_{\ell,k}\}$. Note that, by (6.3), 
$$
-\log (\lambda_{\ell} )\, := 
\,-\,\Sigma_{i=1}^{N_{\ell}}\,\log (\lambda_{\ell, i}) \tau_i \otimes \tau_i^*
$$ 
belongs to $\frak S_{\ell}$. Since the subspace $\frak S_{\ell}$ of $\Bbb R^{N_{\ell}}$ is defined over $\Bbb Q$,
its rational structure  $(\frak S_{\ell})_{\Bbb Q}: = \frak S_{\ell} \cap \Bbb Q^{N_{\ell}}$, 
dense in $\frak S_{\ell}$,
 sits in  $(\frak g_{\ell})_{\Bbb Q} = (\frak g_{\ell})_{\Bbb Z} \otimes_{\Bbb Z} \Bbb Q$,
Take a sequence of points $\hat{\beta}_{\nu} =  \Sigma_{i=1}^{N_{\ell}}
\,\hat{\beta}_{\nu ;i}\tau_i\otimes \tau_i^*  = (\hat{\beta}_{\nu ;1}, \hat{\beta}_{\nu ;2}, 
\dots, \hat{\beta}_{\nu ;N_{\ell}})\in (\frak S_{\ell})_{\Bbb Q}$, 
$\nu = 1,2,\dots$, satisfying $\Sigma_{i=1}^{N_{\ell}}\hat{\beta}_{\nu; i} = 0$ 
such that we have the convergence
$$
\hat{\beta}_{\nu} \to\; -\,\log  (\lambda^{}_{\ell}), 
\quad\text{ as $\nu \to \infty$}.
\leqno{(6.4)}
$$
%\begin{cases}
%&|\exp (\hat{\gamma}_{k;\alpha} ) - \lambda_{\ell,\alpha}| \; \leq \; 
% (\, \min_{1\leq \alpha \leq N_{\ell}}\ell\,\lambda_{\ell,\alpha}\,)_{}^{-1},
%\quad \alpha = 1,2, \dots, N_{\ell},\\
%&
%\leqno{(3.3)}
%$$
%for a sufficiently large $k$, say for $k = k(\ell )$.
Let $a_{\nu}$ be the smallest positive integer  such that $a_{ \nu}\hat{\beta}_{\nu}$ is integral. 
Write $a_{\nu}\hat{\beta}_{\nu}$ as
${\beta}_{\nu} =  
\,({\beta}_{\nu ;1}, {\beta}_{\nu ;2}, \dots, {\beta}_{\nu ;N_{\ell}}) \in\Bbb Z^{N_{\ell}}$ 
for simplicity.
Then ${\beta}_{\nu} $ can be viewed as an element of $(\frak g_{\ell})_{\Bbb Z}$.
We now define 
an algebraic group homomorphism $\psi^{}_{\ell, k} : \Bbb G_m=\{\,t\in \Bbb C^*\} \to \operatorname{SL}(V_{\ell})$ 
by 
$$
\psi^{}_{\ell, \nu} (t)\cdot \tau^{}_{i} \; :=\; t_{}^{-{\beta}^{}_{\nu ;i}}\tau^{}_{i},
\qquad i = 1,2,\dots, N_{\ell}.
$$
Next we identify $X$ with the subvariety $\hat{X}_{\ell} := \hat{\Phi}_{\ell}(X)$ in the projective 
space $\Bbb P^{N_{\ell}-1}(\Bbb C ) = \{(z_1:z_2: \dots :z_{N_{\ell}})\}$ 
by the Kodaira embedding  $\hat{\Phi}_{\ell}$ as in (4.5).
For each $\ell \geq \ell_0$, we observe that $\operatorname{SL}(V_{\ell})$ acts naturally on $W^*_{\ell}$.
%Note that $\psi_{\ell,k}^{} \,=\, \psi_{\ell,k}^{\operatorname{SL}}$.
Then by considering the sequence of test configurations 
$$
\mu^{}_{\ell,\nu} = (\mathcal{X}^{\psi^{}_{\ell, \nu }},
\mathcal{L}^{\psi^{}_{\ell, \nu}},{\psi^{}_{\ell, \nu}}), \qquad \nu = 1,2,\dots,
$$
associated to ${\beta}_{\nu}\in (\frak g_{\ell})_{\Bbb Z}$ (cf.~(3.2)),
and by restricting $\Bbb G_m$ to its real Lie subgroup $\Bbb R_+$,
define a real-valued function 
$\hat{f}^{}_{\ell,\nu} = \hat{f}^{}_{\ell,\nu} (s)$ on $\Bbb R = \{ s \in \Bbb R \}$ by
$$
\hat{f}^{}_{\ell,\nu} (s) \; :=\; 
\begin{cases}
&\delta (\psi^{}_{\ell,\nu} ) \,
\ell^{-n}\log  \| \psi_{\ell,\nu}^{}(t)\cdot \hat{X}_{\ell}\|_{\operatorname{CH}(\rho^{}_{\ell})}, \quad 
\text{ if $\|\psi_{\ell,\nu}\|_{\infty} \neq 0$},\\
&\delta (\psi^{}_{\ell,\nu} ) \,
\ell^{-n}\log  \|  \hat{X}_{\ell}\|_{\operatorname{CH}(\rho^{}_{\ell})}, \;\;\qquad \qquad
\text{ if $\|\psi_{\ell,\nu}\|_{\infty} = 0$},
\end{cases}
$$
where $s\in \Bbb R$ and $t\in \Bbb R_+$ are related by $t = \exp (s/\|\psi_{\ell,\nu}\|_{\infty})$
for  $\|\psi_{\ell,\nu}\|_{\infty} \neq 0$, and we require no relations between $s$ and $t$ 
if  $\|\psi_{\ell,\nu}\|_{\infty} = 0$. Here $\| \psi_{\ell,\nu}\|_{\infty}$ and $\delta (\psi_{\ell,\nu } )$ 
are as in Section 5.
Put %$\dot{f}^{}_{\ell,k}:= df^{}_{\ell,k}/ds$ and 
$\theta_{s;\ell,\nu} := (1/2\pi )\log \{(\Sigma_{i = 1}^{N_{\ell}} \, 
t^{2{\beta}^{}_{\nu ;i}} 
|\hat{\tau}^{}_{i}|^2)^{1/\ell}\}$. 
By identifying $X$ with $\hat{X}_{\ell}$ via $\hat{\Phi}_{\ell}$, we obtain
$$
\psi_{\ell,\nu}^{}(t)^*({\omega}_{\operatorname{FS}}/\ell){}_{|\hat{X}_{\ell}}
= \sqrt{-1}\partial\bar{\partial}\theta_{s;\ell,\nu}, 
\leqno{(6.5)}
$$
where 
${\omega}_{\operatorname{FS}} 
:= (\sqrt{-1}/2\pi ) \partial\bar{\partial}\log \{(\Sigma_{i = 1}^{N_{\ell}}
|z_i|^2)_{}^{1/\ell}\}$.
In view of \cite{Zh} (see also \cite{MT1} and \cite{Sn}), we obtain
$$
d\hat{f}^{}_{\ell,\nu}/ds \; =\; \ell\,\delta (\psi_{\ell,\nu})
\int_X (\partial \theta_{s;\ell,\nu}/\partial s)\,(\sqrt{-1}\partial\bar{\partial}\theta_{s;\ell,\nu })^n.
\leqno{(6.6)}
$$
Put
$\eta^{}_{\ell,\nu} := \|\psi_{\ell,\nu}\|_{\infty}/a_{\nu} 
= \max \{\, |\hat{\beta}_{\nu;i}|/\ell\, ;\, i = 1,2,\dots,N_{\ell}\,\}$, where
%+ a_{\ell, k}^{-1}\inf \|\varphi\|_{\infty}$, 
%where the infimum is taken over all $\ell$-th roots $(\mathcal{Y}, \mathcal{Q}, D, \varphi )$
%of $\mu_{\ell,k}$, and 
for the time being, we vary $\ell$ and $\nu$ independently. Then
$$
(\partial \theta_{s;\ell,\nu}/\partial s )_{\,|s \,=\, -\eta^{}_{\ell,\nu}}\; =\;
\frac{\Sigma_{i =1}^{N_{\ell}}\,\hat{\beta}_{\nu;i} \,
\exp (-2{\hat{\beta}^{}_{\nu ;i}})\, |\hat{\tau}^{}_{i}|^2}
{\pi \ell \,\eta^{}_{\ell,\nu} \,\Sigma_{i =1}^{N_{\ell}}\,
\exp (-2{\hat{\beta}^{}_{\nu ;i}})\,  |\hat{\tau}^{}_{i}|^2}.
\leqno{(6.7)}
$$
Now for each integer $r$, let $O(\ell^{\,r})$ denote a function $u$ 
%or a real $2$-form $\theta$ on $X$ such that 
%$|u |\leq C_0 \ell_j^{\,r}$ or $-C_0 \ell_j^{\,r}\omega
%\leq \theta \leq C_0 \ell_j^{\,r}\omega$, respectively,  
satisfying the inequality $|u| \leq C_0\ell^{\,r}$
for some positive constant $C_0$ independent of %the choice of 
$\nu$, $\ell$, and $i$.  We now fix a positive integer $\ell \gg 1$. Then by (6.4), we obtain
$$
\lambda^{-2}_{\ell,i}\,\exp (-2\hat{\beta}_{\nu ;i})\, -\,1  \; =\;  O(\ell^{-n-2}), \qquad \nu\gg 1.
\leqno{(6.8)}
$$
Moreover, in view of  (6.1) and (6.2), the K\"ahler form
$\omega_{\ell}$ is written as
$(\sqrt{-1}/2\pi )\partial\bar{\partial}\log \{ (\Sigma_{i =1}^{N_{\ell}} 
\lambda^{\,2}_{\ell,i}| \hat{\tau}^{}_{i}|^2)_{}^{1/\ell}\}$. Now by (6.4),
 as $\nu \to \infty$, we have $\sqrt{-1}\partial\bar{\partial}\theta_{s;\ell,\nu\,|\,s\,=\,-\eta^{}_{\ell,\nu}} \to \omega_{\ell}$ in $C^{\infty}$. In particular for $\nu \gg 1$,
$$
\|
\sqrt{-1}\partial\bar{\partial}\theta_{s;\ell,\nu\,|\,s\,=\,-\,\eta^{}_{\ell,\nu}} - \omega_{\ell}\,
\|^{}_{C^5(X)}
\;  = \; O(\ell^{-n-2}).
\leqno{(6.9)}
$$
Hence for each $\ell \geq \ell_0$, we can find a positive integer $\nu (\ell) \gg 1$ such that
both (6.8) and (6.9) hold for $\nu  = \nu (\ell)$. From now on, we assume 
$$
\nu \;=\; \nu (\ell ),
\leqno{(6.10)}
$$
and $\eta^{}_{\ell,\nu }\,=\,\eta^{}_{\ell,\nu (\ell )}$ will be written as $\eta_{\ell}$ for simplicity. 
Then, since $\ell \eta^{}_{\ell} \geq |\hat{\beta}^{}_{\nu ;i}|$ for all $i$, 
we have $(\partial \theta_{s;\ell,\nu }/\partial s )^{}_{|\,s \,=\, -\eta^{}_{\ell}} = O(1)$ by (6.7). Hence 
$$
\int_X (\partial \theta_{s;\ell,\nu}/\partial s)\,\{ \,(\sqrt{-1}\partial\bar{\partial}\theta_{s;\ell,\nu})^n \,-\, \omega_{\ell}^{\,n}\,\}^{}_{\,|\,s\,=\,-\eta^{}_{\ell}} \; = \; O(\ell^{-n-2}).
\leqno{(6.11)}
$$
Put $I_1 := \pi \ell \,\eta^{}_{\ell} \,\Sigma_{i =1}^{N_{\ell}}\,
\lambda^2_{\ell,i}\,  |\hat{\tau}^{}_{i}|^2 $ and $I_2 := \pi \ell \,\eta^{}_{\ell} \,\Sigma_{i =1}^{N_{\ell}}\,
\exp (-2{\hat{\beta}^{}_{\nu ;i}})\, |\hat{\tau}^{}_{i}|^2 $.
Put also $J_1 := \Sigma_{i =1}^{N_{\ell}}\,\hat{\beta}_{\nu ;i}
\lambda^2_{\ell,i} |\hat{\tau}^{}_{i}|^2 $ and $J_2 :=  \Sigma_{i =1}^{N_{\ell}}
\hat{\beta}_{\nu ;i}\exp (-2{\hat{\beta}^{}_{\nu ;i}}) |\hat{\tau}^{}_{i}|^2$.
Then by (6.7), we obtain
$$
\int_X \,(\partial \theta_{s;\ell,\nu}/\partial s )^{}_{\,|s\,=\,-\eta^{}_{\ell}} \;\omega_{\ell}^{\,n}\;
= \; A \,+\, P\, +\,Q,
\leqno{(6.12)}
$$
where $ A :=\int_X  \{(J_2/I_2) - (J_2/I_1)\} \,\omega_{\ell}^{\,n}$,
$P := \int_X \{(J_2/I_1) - (J_1/I_1)\}\,\omega_{\ell}^{\,n}$ and 
$Q := \int_X \,(J_1/I_1)\,\omega_{\ell}^{\,n}$. Note that $J_2/I_2  = O(1)$
by $\ell \eta^{}_{\ell} \geq |\hat{\beta}^{}_{\nu ;i}|$,
while by (6.8), $(I_1 - I_2)/I_1 = O(\ell^{-n-2})$. Then 
$$
A \; =\; \int_X\, \frac{J_2}{I_2}\cdot\frac{I_1 - I_2}{I_1}\;\omega_{\ell}^n\; =\; O(\ell^{-n-2}).
\leqno{(6.13)}
$$
On the other hand by (6.8), $J_2 - J_1 = O(\ell^{-n-2}) (\Sigma_{i =1}^{N_{\ell}}\,
|\hat{\beta}_{\nu ;i}|\,\lambda^2_{\ell,i}\,    |\hat{\tau}^{}_{i}|^2)$.
From this together with $\ell \eta^{}_{\ell} \geq |\hat{\beta}^{}_{\nu ;i}|$,
we obtain 
$$
P \; =\; \int_X \frac{J_2 - J_1}{I_1} \;\omega_{\ell}^{\,n} \; = \; O(\ell^{-n-2}).
\leqno{(6.14)}
$$
By (6.2),  $ I_1 = \pi \ell \,\eta^{}_{\ell} \,\Sigma_{i =1}^{N_{\ell}}\,
  |\hat{\sigma}^{}_{i}|^2$ and $J_1 := \,\Sigma_{i =1}^{N_{\ell}}\,\hat{\beta}_{\nu ;i}
\,  |\hat{\sigma}^{}_{i}|^2 $.  Note also that $a^{}_0 :=  \delta (\psi^{}_{\ell, \nu} )$ satisfies $ 0 < a^{}_0 \leq \ell^n$.
Put $a_1 := c_1(L)^n[X]$.
In view of (6.1) and (6.6), by adding up (6.11), (6.12), (6.13) and (6.14), we see that 
the derivative $d\hat{f}^{}_{\ell,\nu}/ds$ at $s = -\eta^{}_{\ell}$ is written as
$$
\begin{cases}
& %d\hat{f}^{}_{\ell,\nu}/ds\, =\, 
\ell\,a_0\int_X \{\,(\partial \theta_{s;\ell,\nu }/\partial s)\,(\sqrt{-1}\partial\bar{\partial}\theta_{s;\ell,\nu})^n\}^{}_{\,|s\,=\,-\eta^{}_{\ell}}\\
&\;=\;  a^{}_0\,\{\, \ell\,Q  \,+\,O(\ell^{-n-1})\,\}
=\,\int_X \frac{a^{}_0\Sigma_{i =1}^{N_{\ell}}\,\hat{\beta}_{\nu ;i}
\,  |\hat{\sigma}^{}_{i}|_{h_{\ell}}^2\;}{\;\pi \,\eta^{}_{\ell} \,\Sigma_{i =1}^{N_{\ell}}\,
  |\hat{\sigma}^{}_{i}|_{h_{\ell}}^2\;}\, \omega_{\ell}^n \,+\, O(\ell^{-1})\\
  &\; =\; 
 a_0\, a_1 \,(\Sigma_{i =1}^{N_{\ell}}\,\hat{\beta}_{\nu ;i})\,
 (\pi\, \ell \,\eta^{}_{\ell}\, N_{\ell})^{-1} \, +\, O(\ell^{-1}) \; =\;O(\ell^{-1}),  
\end{cases}
\leqno{(6.15)}
$$
where in the last line, we used the equality $\Sigma_{i =1}^{N_{\ell}}\,\hat{\beta}_{\nu ;i} = 0$.
In the next section, the sequence of test configurations $\mu^{}_{\ell,\nu (\ell )} = (\mathcal{X}^{\psi^{}_{\ell, \nu (\ell )}},
\mathcal{L}^{\psi^{}_{\ell, \nu (\ell )}},{\psi^{}_{\ell, \nu (\ell )}})$, $\ell \geq \ell_0$, for $(X,L)$ 
will be considered.

\section{Proof of Main Theorem}

In this section, under the same assumption as in the previous section, we shall show that 
$c_1(L)$ admits an extremal K\"ahler metric. 
Put 
$$
\eta^{}_{\infty}\;  :=\; \sup_{\ell} \; \eta^{}_{\ell},
$$
where the supremum is taken over all positive integers $\ell$ satisfying $\ell \geq \ell_0$.
Then the following cases are possible:

\medskip
 \qquad Case 1:  \; $\eta^{}_{\infty}\, =\, +\infty$.
 \qquad Case 2:  \; $\eta^{}_{\infty} \,<\, +\infty$.
 
 \medskip\noindent
{\em Step \/$1$.}  If Case 1 occurs, then an increasing subsequence $\{\,\ell_j\,; \,j=1,2,\dots \,\}$ 
of  $\{\,\ell\in \Bbb Z\,;\, \ell \geq \ell_0\,\}$ 
can be chosen in such a way that $\{ \eta^{}_{\ell_j}\}$ is a  monotone increasing sequence 
satisfying
$$
\lim_{j\to \infty} \,\eta^{}_{\ell_j}\; =\; +\infty.
\leqno{(7.1)}
$$
For simplicity, the function $\hat{f}^{}_{\ell_j, \nu (\ell_j )}$ will be written as $\hat{f}_{j}$, while  
we write the test configurations 
$$
\mu^{}_{\ell_j, \nu (\ell_j )}= (\mathcal{X}^{\psi^{}_{\ell_j, \nu (\ell_j )}},
\mathcal{L}^{\psi^{}_{\ell_j, \nu (\ell_j )}},{\psi^{}_{\ell_j, \nu (\ell_j )}}), \qquad j=1,2,\dots,
$$
as $\mu_j = (\mathcal{X}_j, \mathcal{L}_j, \psi_j )$.
Now by (6.15), there exists a positive constant $C$ independent of $j$
 such that 
 $$
-\, C /\ell_j\;\leq \;  (d\hat{f}_{j}^{}/ds)\,{}_{| s = -\eta^{}_{\ell_j}}
 $$
for all $j$. On the other hand, for all positive integers $ j' $ satisfying $j' \geq j$, we have $-\eta^{}_{\ell_{j'}}
\leq -\eta^{}_{\ell_j}$ by monotonicity. Since the function $ d\hat{f}_{j'}^{}/ds$ in $s$ is non-decreasing, 
we obtain
$$
 -\, C/\ell_{j'} \;\leq \; d\hat{f}_{j'}^{}/ds\,{}_{|s= -\eta^{}_{\ell_{j'}}}\; \leq \;  d\hat{f}_{j'}^{}/ds\,{}_{|s=-\eta^{}_{\ell_{j}}}.
 \leqno{(7.2)}
$$
We here observe that $ -\, C/\ell_{j'}  \to 0$ as $j' \to \infty$.  It now follows from (7.2) that, for each fixed $j$, 
$$
\varliminf_{j' \to \infty}d\hat{f}_{j'}^{}/ds\,{}_{|s=-\eta^{}_{\ell_{j}}} \; \geq\; 0.
$$
Since the function  $\varliminf_{j' \to \infty}d\hat{f}_{j'}^{}/ds$ in $s$ is non-decreasing, 
we obtain 
$$
\varliminf_{j' \to \infty} d\hat{f}_{j'}^{}/ds \; \geq \; 0
\;\, \,\text{on  \,$\{\,s\in \Bbb R\,;\, s \geq -\eta^{}_{\ell_{j}}\,\}$,}
$$
while this holds for all positive integers $j$. Then by (7.1), 
$\varliminf_{j' \to \infty}d\hat{f}_{j'}^{}/ds$ is, as a function in $s$, nonnegative on the whole 
real line $\Bbb R$.
Hence 
$$
\hat{F}_1( \{\mu_j\}) \; =\; \lim_{s\to -\infty}\{\varliminf_{j' \to \infty}d\hat{f}_{j'}^{}/ds \} \; \geq \; 0.
$$
By Lemma 5.9, ${F}_1( \{\mu_j\}) \geq 0$.
Hence by strong K-stability of $(X,L)$ relative to $T$, we obtain ${F}_1(\{\mu_j \}) = 0$, so that
$\mu_j$ are trivial  for all $j\gg 1$.
This then usually gives us a contradiction. Even if no contradictions occur, from the triviality 
of $\mu_j$ for all $j\gg 1$, we still proceed as follows: 
Since $\psi_j^{}\, =\, \psi^{}_{\ell_j, \nu (\ell_j)}$ are trivial 
for all $j\gg 1$,  by the equality (6.5), we see for all $s\in\Bbb R$ that
$$
\sqrt{-1} \partial\bar{\partial}\theta_{s;\ell_j, \nu (\ell_j )}\; =\;
({\omega_{\operatorname{FS}}}/\ell_j )_{|\hat{X}_{\ell_j}} \; =\; \hat{\Phi}_{\ell_j}^*({\omega_{\operatorname{FS}}}/\ell_j )
\qquad  j \gg 1,
\leqno{(7.3)}
$$
under the identification of $\hat{X}_{\ell}$ with $X$ via $\hat{\Phi}_{\ell}$. 
First by \cite{Ze} (see also \cite{M3}),  
$$
 \|\Phi_{\ell_j}^*({\omega_{\operatorname{FS}}}/\ell_j )\, -\, \omega\,\|^{}_{C^5(X)}\,=\,O(\ell_j^{-2}). 
 \leqno{(7.4)}
$$
In view of the natural isomorphism $\frak t \cong \frak t_{\ell}$, we regard $v_{\ell} \in \frak t_{\ell}$ in (1) of Section 4
as an element of $\frak t$. 
Then we can write $\hat{\Phi}_{\ell_j}^*({\omega_{\operatorname{FS}}}/\ell_j ) - 
\Phi_{\ell_j}^*({\omega_{\operatorname{FS}}}/\ell_j )$ as 
%\begin{align*}
%\hat{\Phi}_{\ell_j}^*({\omega_{\operatorname{FS}}} ) - 
%\Phi_{\ell_j}^*({\omega_{\operatorname{FS}}} ) 
%& = 
$$
\begin{cases}
&(\sqrt{-1}/2\pi )\,\ell_j^{-1}\partial \bar{\partial}\log
\left \{ 1 + 
\frac{\Sigma_{k=1}^{m_j}\Sigma_{\alpha}^{n_{j,k}}(\gamma_{\ell,k} - 1)|\tau_{k,\alpha}|^2}
{\Sigma_{k=1}^{m_j}\Sigma_{\alpha}^{n_{j,k}}|\tau_{k,\alpha}|^2} 
\right \}\\
&= \; (\sqrt{-1}/2\pi )\, \ell_j^{-1} \partial \bar{\partial}(1+ \ell_j^{-1}K_j),
\end{cases}
\leqno{(7.5)}
$$
where $K_j$, $j =1,2,\dots$,  are uniformly bounded Hamiltonian functions 
associated to the holomorphic vector field $\ell^{\,2}_j\, v_{\ell_j}\in \frak t$ 
on the K\"ahler manifold $(X, \Phi_{\ell_j}^*({\omega_{\operatorname{FS}}}/\ell_j ))$,
so that $i_{\ell^{\,2}_j\, v_{\ell_j}}({\omega_{\operatorname{FS}}}/\ell_j ) =
(\sqrt{-1}/2\pi ) \bar{\partial} K_j$. In particular
$$
 \|\hat{\Phi}_{\ell_j}^*({\omega_{\operatorname{FS}}}/\ell_j )\, -\,\Phi_{\ell_j}^*({\omega_{\operatorname{FS}}}/\ell_j )\,\|^{}_{C^5(X)}\,=\,O(\ell_j^{-2}). 
 \leqno{(7.7)}
$$
Then by (7.4) and (7.6), $\| \hat{\Phi}_{\ell_j}^*({\omega_{\operatorname{FS}}}/\ell_j )\,
 -\, \omega\,\|^{}_{C^5(X)}\,=\,O(\ell_j^{-2})$.
Hence by applying (7.3) to $s = -\eta^{}_{\ell_j}$, we see from (6.9) that
$$
\|\,\omega \,-\,\omega^{}_{\ell_j}\,\|^{}_{C^5(X)} \; =\; O(\ell_j^{-2}), \qquad j \gg 1.
\leqno{(7.8)}
$$
Let $S_{\omega}$ be the scalar curvature function for $\omega$. Then by \cite{Lu} (see also \cite{Ze}),
we obtain the following asymptotic expansion:
$$
1+ (S_{\omega}/2)\ell_j^{-1} + O(\ell_j^{-2})\; =\; 
\Sigma_{\alpha =1}^{N_{\ell_j}}
(n!/\ell_j^n ) \,|\tau^{}_{\ell_j, \alpha}|^2_{h}\; =\; (n!/\ell_j^n ) \,B_{\ell_j}(\omega ).
\leqno{(7.9)}
$$
On the other hand, for $\ell \gg 1$,  we observe that $N_{\ell}$ is a polynomial in $\ell$. Since each $\omega_{\ell_j}$ is polybalanced, by setting $\ell = \ell_j$ in (4) of Section 4, we obtain 
$$
1 +  n! \hat{H}_{j}  \ell^{-1}_j + O( \ell_j^{-2}) = \Sigma_{\alpha =1}^{N_{\ell_j}}
(n!/\ell_j^n ) |\sigma^{}_{\ell_j, \alpha}|^2_{h^{}_{\ell_j}}
 = (n!/\ell_j^n ) B_{\ell_j}(\omega_{\ell_j} ),
\leqno{(7.10)}
$$
where $\hat{H}_{\ell_j} = H_{\ell_j} + C^{}_0$ for some constant $C_0$ independent of $j$. 
In view of (7.6), 
by comparing (7.7) with (7.8), we conclude from the uniformly boundedness of $H_{\ell_j}$, $j =1,2,\dots$, 
we see 
that $S_{\omega}$ is an Hamiltonian function for $(X,\omega)$. Hence $\omega$ is an extremal K\"ahler metric in the class $c_1(L)_{\Bbb R}$.

\medskip\noindent
{\em Step \/$2$.} Suppose that Case 2 occurs. Put ${\lambda}'_{\ell, i} := -(1/\ell )\log \lambda_{\ell, i}$. 
Then by (6.4), we may assume that $\nu = \nu (\ell )$ in (6.10) is chosen in such a way that 
$$
\hat{\beta}_{\nu (\ell );i} -1\, \leq \; \ell\, {\lambda}'_{\ell, i}\; \leq \,\hat{\beta}_{\nu (\ell );i} +1,
\qquad i = 1,2,\dots, N_{\ell},
\leqno{(7.11)}
$$
for all $\ell$ with $\ell \geq \ell_0$. 
%Let $i$ be arbitrary integers satisfying
%$1\leq i  \leq N_{\ell}$.
By (7.11) together with the definition of $\eta_{\ell,\nu (\ell ) }$ in Section 6,
we see from the inequality $\nu_{\infty}  < +\infty$ that
$$
|{\lambda}'_{\ell, i}| \;\leq\; C_7, \qquad i=1,2,\dots,N_{\ell},
$$
where $C_7$ is a positive real constant independent of the choice of $\ell$ and $i$.
Then by taking an $\ell$-th root (cf.~\cite{M1}, Section 5) for $\mu_{\ell,\nu (\ell )}$ in Section 6,
and by the same argument as in Step 2 of \cite{M1}, Section 4,
we have 
$$
\omega_{\ell_j} \to \omega_{\infty},
\qquad\text{as $j \to \infty$},
$$
for some K\"ahler form $\omega_{\infty}$ on $X$.
Then it is easy to see that
$\omega_{\infty}$ is an extremal K\"ahler metric in the class $c_1(L)_{\Bbb R}$,
as required.

\section{Appendix}

For a polarized algebraic manifold $(X,L)$, let $T$ be as in the introduction.
Then the purpose of this appendix is to show the following:

\medskip\noindent
{\bf Theorem 8.1.} 
{\em If $c_1(L)_{\Bbb R}$ admits an extremal K\"ahler metric, then 
$(X,L)$ 
is strongly K-semistable relative to $T$.} %relative to $G$.} 

\medskip
In \cite{M2}, we shall show a stronger result that
the polarized algebraic manifold $(X,L)$ in Theorem B is actually strongly K-stable relative to $T$.

\medskip
In this appendix, we use the same notation as in previous sections, 
where we assume that $\omega = c_1(L;h)$ is an extremal  K\"ahler metric in the class $c_1(L)_{\Bbb R}$
on $X$. 
In view of \cite{M2} and \cite{MN1}, $(X,L)$ is asymptotically Chow-stable relative to $T$ (see also \cite{ST} and their subsequent works).
For $\omega$ and $h$ above, we consider the Hermitian metric $\rho_{\ell}$ on $V_{\ell}$ 
as in (2.1).
Then the scalar curvature function
$S_{\omega}$ for the K\"ahler form $\omega$
% \; =\;  n\, c_1(M)c_1(L)^{n-1}[X]/c_1(L)^n[X]
is a Hamiltonian function for an extremal K\"ahler vector field $\mathcal{V}\in H^0(X, \mathcal{O}(T_X))$
in the sense that
$i^{}_{\mathcal{V}}\,\omega \, =\, (\sqrt{-1}/2\pi ) \bar{\partial} S_{\omega}$.
Note that the average $\bar{S}_{\omega}$  of the scalar curvature $S_{\omega}$ on 
the K\"ahler manifold $(X, \omega )$ 
is 
$n\, c_1(M)c_1(L)^{n-1}[X]/c_1(L)^n[X]$. 
Let 
$$
\mu_j = (\mathcal{X}_j, \mathcal{L}_j, \psi_j ), \qquad j=1,2,\dots, 
$$
be a sequence of test configurations 
for $(X,L)$ as in (2.2), so that the exponent $\ell_j$ for $\mu_j$ satisfies $\ell_j \to +\infty$ as $j \to \infty$.
Put $\hat{S}_{\omega} := \bar{S}_{\omega}+c_1(L)^n[X]^{-1}2\pi F(\mathcal{V})$, where 
$F$ is the classical Futaki character.
Let $\hat{B}_j (\omega )$  be the as in (5.3).
Then by \cite{MT1}, \cite{MT'}, \cite{MT2} and \cite{ST},
$\hat{B}_j (\omega )$  is written as
$$
\hat{B}_j (\omega )\; =\; (1/n!)\, \{ \ell_j^n \,+\,(\hat{S}_{\omega}/2 ) \ell_j^{n-1} \,+\, R_j \ell_j^{n-2}\},
\leqno{(8.2)}
$$
where $R_j$ is a function on $X$ satisfying $\|R_j\|^{}_{C^2(X)}\leq C$. 
%Now for each integer $r$, 
%we denote by $O(\ell_j^{\,r})$ a real function $\varphi$ 
%or a real $2$-form $\theta$ on $X$ such that 
%$|\varphi |^{}\leq C' \ell_j^{\,r}$ or $-C' \ell_j^{\,r}\omega
%\leq \theta \leq C' \ell_j^{\,r}\omega$, respectively,  
%for some positive real constant $C'$ 
%independent of $j$. 
Then by taking 
$(\sqrt{-1}/2\pi) \partial\bar{\partial}\log$ of both sides of (8.2), we obtain
$$
\hat{\Phi}_j^*\omega_{\operatorname{FS}} - \ell_j \omega\; =\; O(\ell_j^{-2}).
\leqno{(8.3)}
$$
Hence by (8.1) and (8.2), we see from (5.6) that
\begin{align*}
&d\hat{f}_{j}/ds \, {}_{|s=0} = (\|\psi_j\|_1^{}\ell_j^n )^{-1}\int_X
\frac{n!\,\Sigma_{k=1}^{m_j}\Sigma_{\alpha =1}^{n_{j,k}}\, b^{}_{j;k,\alpha} |\hat{\tau}_{k,\alpha}|_h^2}{   \ell_j^n +(\hat{S}_{\omega}/2 ) \ell_j^{n-1} + O(\ell_j^{n-2})}  \,\{\ell_j\omega + O(\ell_j^{-2})\}^{\,n}\\
%&=\; (|b|\ell_j^{\,n})^{-1}\int_X
%\frac{n!\,\Sigma_{\alpha =1}^{N_j}\, b_{\alpha} |{\sigma}_{\alpha}|_h^2}
%{   \ell_j^n \,+\,(S_{\omega}/2 ) \ell_j^{n-1} + O(\ell_j^{n-2})}  
%\,\{\ell_j\omega + O(\ell_j^{-2})\}^{\,n}\\
&=\; \frac{n!\,\ell_j \|b_j\|_1^{-1}}{\ell_j^n \,+\,(\hat{S}_{\omega}/2 ) \ell_j^{n-1} }
\int_X \frac{\Sigma_{k=1}^{m_j}\Sigma_{\alpha =1}^{n_{j,k}}\, b^{}_{j;k,\alpha} |\hat{\tau}_{k,\alpha}|_h^2}{  1 \,+\,O(\ell_j^{-2})} 
\,\{\ell_j\omega + O(\ell_j^{-2})\}^{\,n}.
\end{align*}
Since $\int_X \Sigma_{k=1}^{m_j}\Sigma_{\alpha =1}^{n_{j,k}}\, b^{}_{j;k,\alpha} |\hat{\tau}_{k,\alpha}|_h^2\, \omega^n = 
\Sigma_{k=1}^{m_j}\Sigma_{\alpha =1}^{n_{j,k}}\, b^{}_{j;k,\alpha}\gamma_{j,k}= 0$ as in the proof of 
Proposition 5.2,  there exists
a positive constant $C''$ independent of $j$ such that, for all $j \gg 1$, we obtain
\begin{align*}
|d\hat{f}_{j}/ds| \, {}_{|s=0} \; &\leq \; C''\,\ell_j \|b_j\|_1^{-1}\int_X (\Sigma_{k=1}^{m_j}\Sigma_{\alpha =1}^{n_{j,k}}\, |b^{}_{j;k,\alpha}|\, |\hat{\tau}_{k,\alpha}|_h^2)\,O(\ell_j^{-2})\,\omega^n\\
&=\; C''\,\ell_j \|b_j\|_1^{-1}\,(\Sigma_{k=1}^{m_j}\Sigma_{\alpha =1}^{n_{j,k}}\, |b^{}_{j;k,\alpha}|)\,O(\ell_j^{-2})
\; =\; O(\ell_j^{-1}).
\end{align*}
%where the last equality follows from Remark 4.6. 
Then by $\lim_{j\to \infty}\ell_j^{-1} =0$,
 we obtain $\lim_{j\to\infty} d\hat{f}_{j}/ds \,{}_{|s=0} 
=0$. 
%In particular, $\eta$ in Section 2 is nonpositive.
Since the function $\varliminf_{j\to \infty}  d\hat{f}_{j}/ds$ in $s$ is non-decreasing, 
 we see 
for all $\{\mu_j\}\in\mathcal{M}$ that
$$
F_1(\{\mu_j\}) \; =\; \lim_{s\to -\infty} \{
\varliminf_{j\to \infty}  d\hat{f}_{j}/ds\}
\; 
%\; \leq \; F_{1,h}(\{\mu_j\}) 
 \leq \; \varliminf_{j\to \infty} d\hat{f}_{j}/ds\,{}_{|s=0}
\; =\; 0,
$$
as required. This completes the proof of Theorem 8.1.

%%%%%%%%%%%%%%%%%%%%%%%%%%%%%
\bigskip\noindent
{\footnotesize
{\sc Department of Mathematics}\newline
{\sc Osaka University} \newline
{\sc Toyonaka, Osaka, 560-0043}\newline
{\sc Japan}}
%%%%%%%%%%%%%%%%%%%%%%%%%%%%%

\end{document}